\documentclass[a4paper,10pt]{amsart}
\usepackage[latin1]{inputenc}
\usepackage{amssymb}
\usepackage{mathrsfs}
\usepackage[all]{xy}
\usepackage{verbatim}

\title{Strong homotopy properads}
\author{Johan Granåker}
\newcommand{\cC}{\mathcal{C}}
\newcommand{\cE}{\mathcal{E}}
\newcommand{\cF}{\mathcal{F}}
\newcommand{\cG}{\mathcal{G}}
\newcommand{\cP}{\mathcal{P}}
\newcommand{\ra}{\rightarrow}
\newcommand{\la}{\leftarrow}
\newcommand{\Ra}{\Rightarrow}
\newcommand{\xra}{\xrightarrow}
\newcommand{\xla}{\xleftarrow}

\DeclareMathOperator{\Id}{Id}
\DeclareMathOperator{\Bij}{Bij}
\DeclareMathOperator{\coH}{H}
\newcommand{\iso}{\cong}
\newcommand{\decor}[2]{#1\langle#2\rangle}
\theoremstyle{plain}
\newtheorem{thm}{Theorem}
\newtheorem{proposition}{Proposition}
\newtheorem{lemma}{Lemma}
\theoremstyle{definition}
\newtheorem*{definition}{Definition}
\newtheorem*{example}{Example}
\newcommand{\bde}{\begin{definition}}
\newcommand{\ede}{\end{definition}}
\newcommand{\bth}{\begin{thm}}
\renewcommand{\eth}{\end{thm}}
\newcommand{\bpr}{\begin{proposition}}
\newcommand{\epr}{\end{proposition}}
\newcommand{\ble}{\begin{lemma}}
\newcommand{\ele}{\end{lemma}}
\newcommand{\bpro}{\begin{proof}}
\newcommand{\epro}{\end{proof}}
\newcommand{\bex}{\begin{example}}
\newcommand{\eex}{\end{example}}
\begin{document}

\begin{abstract}
In this paper, we define the notion of strong homotopy properads and prove that this structure transfers over left homotopy inverses. We give explicit formulae for the induced structure.
\end{abstract}

\maketitle
\section{Introduction}\hspace{1pt}

The notion of properad was introduced by Vallette in \cite{Vallette2004}. The word properad is a contraction of the words prop and operad, hinting at its place `in between' these two structures. While operads are based on rooted trees and props on arbitrary graphs with flow, properads are based on connected graphs with flow. Hence, as for operads there is only the so-called vertical composition, corresponding to grafting output legs of one graph to input legs of another. As for associative algebras and operads this composition is strictly associative (in the generalized sense of compositions of graphs). For algebras, relaxing the associativity constraint gives the notion of $A_\infty$-algebra, or strongly homotopy associative algebra, originally considered by Stasheff \cite{Stasheff1963}. An $A_\infty$-algebra comes equipped with a coherent hierarchy of `higher homotopies', meaning a family of morphisms each of which is a homotopy for the associativity of the previous. One shows that giving a vector space $V$ the structure of an $A_\infty$-algebra is the same as giving a codifferential on the cofree coalgebra on $V$ shifted, i.e.~a degree $1$, square zero coderivation on the tensor algebra of $V[1]$. In this paper we study an analogous graph construction for properads. A properad is a $\Sigma$-bimodule with a strictly associative composition product, and we want to relax this in the same way as for algebras. Hence, we define a structure of strong homotopy properad on a $\Sigma$-bimodule $\cE$ as a degree $1$, square zero coderivation on the cofree coproperad on $\cE$ shifted. We then interpret this structure in terms of graphs. The notion of strong homotopy operad was considered in the thesis of van der Laan \cite{Laan2004a}, see \cite{Laan2002} and \cite{Laan2004b}.

Let $E$ be a differential graded submodule of a differential graded associative algebra $A$, not necessarily closed under multiplication. It is well-known that if there is an isomorphism of cohomology groups, $\coH(E)\iso\coH(A)$, then there is a structure of $A_\infty$-algebra on $E$ induced by the algebra structure on $A$, see \cite{Kadeishvili}, \cite{Kontsevich2001}, \cite{Markl2004}. In \cite{Merkulov1999} explicit formulae for the defining homotopies are given. Moreover, the natural inclusion of $E$ into $A$ extends to a morphism of $A_\infty$-algebras. More generally, if $g:A\ra E$ and $f:E\ra A$ are such that the composition $fg$ is homotopic to the identity on $A$, then there is an induced $A_\infty$-structure on $E$, such that $f$ extends to a morphism of $A_\infty$-algebras. The main purpose of our paper is to generalize these results to properads, and since they are special cases, the results apply also to operads. Let $\cP$ be a properad and $\cE$ a $\Sigma$-bimodule, and let $g:\cP\ra\cE$ and $f:\cE\ra\cP$ be such that $fg$ is homotopic to the identity on $\cP$. In this situation we acquire \emph{explicit} formulae for an induced strong homotopy properad structure on $\cE$ such that, in the special case when $\cP(m,n)$ is zero except for $m=n=1$, our formulae reduce precisely to the ones in \cite{Merkulov1999}. We also get formulae for the extension of $f$ to a morphism of strong homotopy properads. Finally we show that these results remain valid if $\cP$ is replaced by an sh properad.

All modules are assumed to be differential graded, unless otherwise stated, and over a field $\mathbb{K}$. We work with the cohomological convention, that is, all differentials are of degree $+1$. Tensor products are always over $\mathbb{K}$. The cardinality of the set $X$ is denoted $|X|$, and the degree of the homogeneous element $a$ is denoted $|a|$. As usual, it will be clear from the context what is meant. Everywhere we apply Koszul sign rules, meaning when $a$ is moved passed $b$ the sign $(-1)^{|a||b|}$ appears.

In Section $2$ we present the necessary background for stating and proving the result of the paper. In $2.1$ we describe the graphs upon which our presentation of properads, and of strong homotopy properads, is based. Only a certain level of structure of these graphs are relevant to us, and we describe our scheme for picturing graphs in $2.2$. Sections $2.3$ and $2.4$ concerns the process of contracting graphs, to which the composition of a properad is related. In $2.5$ we recall the notion of a $\Sigma$-bimodule, and of graphs decorated with elements of such a module. Sections $2.6$ and $2.7$ contains the definition of a properad and coproperad, and of the bar construction. In Section $3.1$ we give the definition of a strong homotopy properad, and in $3.2$ we go on to prove the main statement by explicitly constructing the induced codifferential. In $3.3$ we construct the above mentioned extension of the morphism $f$, and in $3.4$ we show how the proofs of the results of the previous two sections can be modified to the situation when $\cP$ is an sh properad.

The author would like to express his gratitude towards his advisor S.~A.~Merkulov. We are also much indebted to the referees for many useful comments.

\section{Graphs, properads and bar construction}

\subsection{Graphs}\hspace{1pt}

\bde
A \emph{graph} $G$ is a triple $(\Gamma,\sigma,\lambda)$, where $\Gamma$ is a finite set, $\sigma$ an involution and $\lambda$ a partition. Elements of $\Gamma$ are called \emph{flags}, or \emph{half-edges}.
\ede

The \emph{vertices} of a graph $G$, $\mathbf{v}(G)$, is the set of blocks of $\lambda$. The vertex of a flag $\gamma$, $\mathbf{v}(\gamma)$, is the block of which $\gamma$ is an element. The \emph{edges} of $G$, $\mathbf{e}(G)$, is the set of pairs of flags forming two-element orbits of $\sigma$. The legs of the graph, $\mathbf{l}(G)$, is the set of fixed points for $\sigma$.

\bex
Let $G$ be the graph $(\Gamma,\sigma,\lambda)$, where
\begin{align*}
\Gamma=&\{a,b,c,d,e,f,g,h,i,j,k,l\}, \\
\sigma=&(be)(cf)(di)(gh)(kl), \\
\lambda=&\{\{a,b,c,d\},\{e,f,g\},\{h,i,j,k,l\}\}.
\end{align*}
Then the following is a picture of $G$.
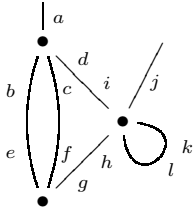
$$
\xygraph{
!{<0pt,0pt>;<30pt,0pt>:<0pt,30pt>::},
!{(0,-1)}*+{\bullet}="a",
!{(0,1)}*+{\bullet}="b",
!{(1,0)}*+{\bullet}="c",
!{(0,1.5)}="d",
!{(1.5,1)}="e",
!{(0.2,1.3)}*+{\scriptstyle{a}},
!{(-0.4,0.4)}*+{\scriptstyle{b}},
!{(-0.4,-0.4)}*+{\scriptstyle{e}},
!{(0.3,-0.4)}*+{\scriptstyle{f}},
!{(0.3,0.4)}*+{\scriptstyle{c}},
!{(0.5,0.8)}*+{\scriptstyle{d}},
!{(0.8,0.5)}*+{\scriptstyle{i}},
!{(1.4,0.5)}*+{\scriptstyle{j}},
!{(0.5,-0.8)}*+{\scriptstyle{g}},
!{(0.8,-0.5)}*+{\scriptstyle{h}},
!{(1.8,-0.3)}*+{\scriptstyle{k}},
!{(1.6,-0.6)}*+{\scriptstyle{l}},
"a"-@/_/"b" "a"-@/^/"b" "a"-"c" "c"-"b" "b"-"d" "c"-"e" "c"-@(d,r)"c"
}
$$
\eex

A morphism of graphs, $G\ra G'$, is a composition of relabelings of flags and contractions of edges. The edge $(ij)$ in $G$ determines a morphism $G\ra G/(ij)$, the contraction of $(ij)$, where, if $G=(\Gamma,\sigma'(ij),\lambda)$,
$$
G/(ij)=(\Gamma\smallsetminus\{i,j\},\sigma',\lambda\smallsetminus(\mathbf{v}(i)\cup\mathbf{v}(j))\cup\{(\mathbf{v}(i)\smallsetminus\{i\})\cup(\mathbf{v}(j)\smallsetminus\{j\})\}).
$$

\bex
The contraction of the edge $(gh)$ in the graph above is the morphism
$$
G=
\xygraph{
!{<0pt,0pt>;<30pt,0pt>:<0pt,30pt>::},
!{(0,-1)}*+{\bullet}="a",
!{(0,1)}*+{\bullet}="b",
!{(1,0)}*+{\bullet}="c",
!{(0,1.5)}="d",
!{(1.5,1)}="e",
!{(0.2,1.3)}*+{\scriptstyle{a}},
!{(-0.4,0.4)}*+{\scriptstyle{b}},
!{(-0.4,-0.4)}*+{\scriptstyle{e}},
!{(0.3,-0.4)}*+{\scriptstyle{f}},
!{(0.3,0.4)}*+{\scriptstyle{c}},
!{(0.5,0.8)}*+{\scriptstyle{d}},
!{(0.8,0.5)}*+{\scriptstyle{i}},
!{(1.4,0.5)}*+{\scriptstyle{j}},
!{(0.5,-0.8)}*+{\scriptstyle{g}},
!{(0.8,-0.5)}*+{\scriptstyle{h}},
!{(1.8,-0.3)}*+{\scriptstyle{k}},
!{(1.6,-0.6)}*+{\scriptstyle{l}},
"a"-@/_/"b" "a"-@/^/"b" "a"-"c" "c"-"b" "b"-"d" "c"-"e" "c"-@(d,r)"c"
}
\ra
\xygraph{
!{<0pt,0pt>;<30pt,0pt>:<0pt,30pt>::},
!{(0,-1)}*+{\bullet}="a",
!{(0,1)}*+{\bullet}="b",
!{(1,0)}="c",
!{(0,1.5)}="d",
!{(1.5,1)}="e",
!{(0.2,1.3)}*+{\scriptstyle{a}},
!{(-0.6,0.4)}*+{\scriptstyle{b}},
!{(-0.6,-0.4)}*+{\scriptstyle{e}},
!{(0.3,-0.4)}*+{\scriptstyle{i}},
!{(0.3,0.4)}*+{\scriptstyle{d}},
!{(-0.1,0.4)}*+{\scriptstyle{c}},
!{(-0.1,-0.4)}*+{\scriptstyle{f}},
!{(0.8,-0.5)}*+{\scriptstyle{j}},
!{(0.7,-1.3)}*+{\scriptstyle{k}},
!{(-0.1,-1.5)}*+{\scriptstyle{l}},
"a"-@/_15pt/"b" "a"-@/^15pt/"b" "a"-"b" "a"-"c" "b"-"d" "a"-@(d,r)"a"
}
=G/(gh).
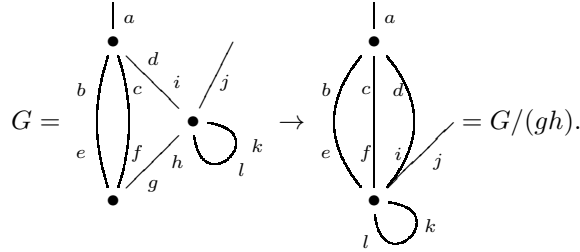
$$
\eex

\bde
We denote by $\cG\hspace{-1pt}\mathit{raph}_k$ the set of isomorphism classes of graphs with $k$ vertices. Graphs with only one vertex are called \emph{corollas}.
\ede

\bde
Two vertices, $v$ and $v'$, are \emph{connected} if there is a sequence of edges,
$$
(\{\gamma_0,\gamma'_0\},\ldots,\{\gamma_n,\gamma'_n\}),
$$
such that $\mathbf{v}(\gamma_0)=v$, $\mathbf{v}(\gamma'_i)=\mathbf{v}(\gamma_{i+1})$ for $0\leq i<n$ and $\mathbf{v}(\gamma'_n)=v'$. We call such a sequence a \emph{chain} between $v$ and $v'$. A graph is connected if all pairs of vertices are connected.
\ede

\bde
A graph $G$ is \emph{directed} if a partition of each vertex into a disjoint union $v=in(v)\sqcup out(v)$ is chosen such that for each edge $\{\gamma,\gamma'\}$ we have $\gamma\in out(\mathbf{v}(\gamma))$ if and only if $\gamma'\in in(\mathbf{v}(\gamma'))$. Hence, we have an ordering of each edge $(\gamma^{out},\gamma^{in})$, in pictures this is drawn as an arrow from $\mathbf{v}(\gamma^{out})$ to $\mathbf{v}(\gamma^{in})$.
\ede

The fact that $G$ is directed implies that the set of legs is partitioned into a disjoint union $\mathbf{l}(G)=out(G)\sqcup in(G)$, where
\begin{align*}
out(G)=&\{\gamma\in\mathbf{l}(G);\gamma\in out(\mathbf{v}(\gamma))\}, \\
in(G)=&\{\gamma\in\mathbf{l}(G);\gamma\in in(\mathbf{v}(\gamma))\}.
\end{align*}

\bde
We denote by $\cG\hspace{-1pt}\mathit{raph}(m,n)$ the set of isomorphism classes of graphs $G$ with $|out(G)|=m$ and $|in(G)|=n$. If $G\in\cG\hspace{-1pt}\mathit{raph}(m,n)$, we say that $G$ is an \emph{$(m,n)$-graph}.
\ede

\bde
Given a pair of vertices $v$ and $v'$, we call the set of all directed edges from $v$ to $v'$ a \emph{thick edge}. We denote the set of thick edges by $\mathbf{te}(G)$, and depict a thick edge $\epsilon$ from $v$ to $v'$ as $v\xra{\epsilon}v'$.
\ede

\bde
Given two vertices, $v$ and $v'$, a \emph{directed path} from $v$ to $v'$ is a sequence of thick edges
$$
v=v_0\xra{\epsilon_1}v_1\xra{\epsilon_2}\cdots\xra{\epsilon_n}v_n=v'.
$$
A \emph{directed cycle} at a vertex $v$ is a non-empty directed path from $v$ to $v$.
\ede

\bde
A \emph{subgraph} $H$ of a graph $G=(\Gamma,\sigma,\lambda)$ is a triple $(\Gamma'\subset\Gamma,\sigma',\lambda')$ such that the vertices of $H$ is a subset of the vertices of $G$, and
\begin{enumerate}
\item if $\gamma,\sigma\gamma\in\Gamma'$, then $\sigma'\gamma=\sigma\gamma$,
\item if $\gamma\in\Gamma',\sigma\gamma\not\in\Gamma'$, then $\sigma'\gamma=\gamma$.
\end{enumerate}
\ede

\bde
A \emph{tree} is a graph such that for any pair of vertices, $v$ and $v'$, there is a unique chain from $v$ to $v'$.
\ede

Choosing a leg of a tree, and a direction of that leg, induces a direction of all edges `by gravity'. The distinguished leg $\gamma_r$ is called the \emph{root}, and the rest of the legs are called \emph{leaves}. We usually choose the direction \emph{towards} the root, then the number of incoming edges to a vertex is called its \emph{arity}.

\bde
A tree $t$ is called \emph{binary} if each vertex consists of precisely three flags. For a binary tree
$$
t=
\xygraph{
!{<0pt,0pt>;<10pt,10pt>:<-10pt,10pt>::}
!{(-0.5,0.5)}*+{\scriptstyle{\epsilon_1}},
!{(0.5,-0,5)}*+{\scriptstyle{\epsilon_2}},
!{(1,2)}*+{\scriptstyle{t_1}},
!{(2,1)}*+{\scriptstyle{t_2}},
!{(-0.2,-0.7)}*+{\scriptstyle{\gamma_r}},
!{(0,0)}-!{(1,0)},
!{(1,0)}-@{.}!{(2,0.5)},
!{(1,0)}-@{.}!{(1.5,1)},
!{(0,0)}-!{(0,1)},
!{(0,1)}-@{.}!{(0.5,2)},
!{(0,1)}-@{.}!{(1,1.5)},
!{(0,0)}-!{(-0.5,-0.5)}
},
$$
we will denote by $\epsilon_1$ and $\epsilon_2$ the first two internal edges as in the picture. We call the subtrees $t_1$ and $t_2$, with roots $\epsilon_1$ and $\epsilon_2$, the \emph{principal subtrees} of $t$.
\ede

\bde
Let $G_1=(\Gamma_1,\sigma_1,\lambda_1)$ and $G_2=(\Gamma_2,\sigma_2,\lambda_2)$ be two graphs. Let $i_1,\ldots,i_n\in in(G_1)$ and $o_1,\ldots,o_n\in out(G_2)$ be legs. We denote by $G_1\circ_{o_1,\ldots,o_n}^{i_1,\ldots,i_n} G_2$ the graph whose set of flags is $\Gamma_1\sqcup\Gamma_2$, with involution
$$
\sigma(\gamma)=
\begin{cases}
\sigma_1(\gamma) & \text{if }\gamma\in\Gamma_1\smallsetminus\{i_1,\ldots,i_n\}, \\
\sigma_2(\gamma) & \text{if }\gamma\in\Gamma_2\smallsetminus\{o_1,\ldots,o_n\}, \\
o_j & \text{if }\gamma=i_j, \\
i_j & \text{if }\gamma=o_j,
\end{cases}
$$
and with partition $\lambda_1\cup\lambda_2$. We say that the leg $o_j$ has been \emph{grafted} to the leg $i_j$.
\ede

In the situation when we have a graph $G$, split into two subgraphs $G_1$ and $G_2$ such that $G$ can be considered as the result of grafting $G_1$ and $G_2$ along some set of legs, we write, for short, $G=G_1\circ G_2$.

\bde
Let $G$ be an $(m,n)$-graph. We say that $G$ is \emph{labeled} if we have bijections
\begin{align*}
h_{out}(G):out(G)\ra&\{1,\ldots,m\}, \\
h_{in}(G):in(G)\ra&\{1,\ldots,n\},
\end{align*}
numbering the output- and input-legs of the graph.
\ede

\bde
We denote by $\cG$ the set of isomorphism classes of directed, connected, labeled graphs without directed cycles. The meanings of $\cG_k$ and $\cG(m,n)$ are the same as for $\cG\hspace{-1pt}\mathit{raph}_k$ and $\cG\hspace{-1pt}\mathit{raph}(m,n)$. To $\cG$ we also add the \emph{trivial} graph consisting of only a directed edge considered to be in $\cG(1,1)$. It will play the role of unit element in the free properad.
\ede

From now on we only consider graphs in $\cG$.

\subsection{Pictorial representation of graphs}\hspace{1pt}

When depicting graphs we draw only thick edges. Also, we do not draw legs. Hence, many different graphs will have the same pictorial representation. For example, both of the graphs
$$
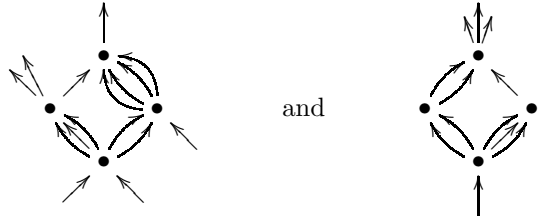
\xygraph{
!{<0pt,0pt>;<10pt,10pt>:<-10pt,10pt>::},
!{(-1,-1)}*+{\bullet}="1",
!{(-1,1)}*+{\bullet}="2",
!{(1,-1)}*+{\bullet}="3",
!{(1,1)}*+{\bullet}="4",
!{(-2.5,-1)}="5",
!{(-1,-2.5)}="6",
!{(-1,2.5)}="7",
!{(-0.5,2.5)}="8",
!{(1,-2.5)}="9",
!{(2,2)}="10",
"1":@/_0.3pc/"2" "1":"2" "1":@/^0.3pc/"2"
"1":@/_0.3pc/"3" "1":@/^0.3pc/"3"
"2":"4"
"3":@/_0.7pc/"4" "3":@/_0.3pc/"4" "3":@/^0.3pc/"4" "3":@/^0.7pc/"4"
"5":"1" "6":"1"
"2":"7" "2":"8"
"9":"3"
"4":"10"
}
\hspace{1cm}\text{ and }\hspace{1cm}
\xygraph{
!{<0pt,0pt>;<10pt,10pt>:<-10pt,10pt>::},
!{(-1,-1)}*+{\bullet}="1",
!{(-1,1)}*+{\bullet}="2",
!{(1,-1)}*+{\bullet}="3",
!{(1,1)}*+{\bullet}="4",
!{(-2,-2)}="6",
!{(2,2)}="7",
!{(1.5,2)}="8",
!{(2,1.5)}="9",
"1":@/_0.3pc/"2" "1":@/^0.3pc/"2"
"1":@/_0.3pc/"3" "1":"3" "1":@/^0.3pc/"3"
"2":@/_0.3pc/"4" "2":@/^0.3pc/"4"
"3":"4"
"6":"1"
"4":"7" "4":"8" "4":"9"
}
$$
would be represented as
$$
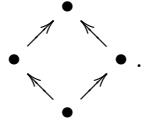
\xygraph{
!{<0pt,0pt>;<10pt,10pt>:<-10pt,10pt>::}
!{(-1,-1)}*+{\bullet}="1",
!{(-1,1)}*+{\bullet}="2",
!{(1,-1)}*+{\bullet}="3",
!{(1,1)}*+{\bullet}="4",
"1":"2" "1":"3" "2":"4" "3":"4"
}.
$$

Only this level of structure of the graphs is relevant to our constructions.

\subsection{Contraction sequences for a graph, and associated trees}\hspace{1pt}

\bde
Given a graph $G=(\Gamma,\lambda,\sigma)$, and a thick edge $v_i\xla{\epsilon}v_j$, we denote by $G/\epsilon$ the graph with the vertices $v_i,v_j$ of $G$ replaced by the single vertex $v_i\sqcup v_j\smallsetminus\{\text{flags in }\epsilon\}$. We say that $G$ is \emph{contracted} along $\epsilon$.
\ede

We commonly identify the vertices and edges of $G/\epsilon$ not affected by the contraction with their pre-images in $G$.

\bde
A thick edge, connecting a pair of vertices $v_i\xla{\epsilon}v_j$, is said to be \emph{admissible} if there are no directed cycles in $G/\epsilon$.
\ede

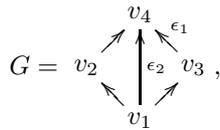
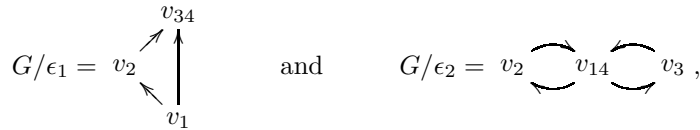
\bex
If
$$
G=
\xygraph{
!{<0pt,0pt>;<10pt,10pt>:<-10pt,10pt>::}
!{(-1,-1)}*+{v_1}="1",
!{(-1,1)}*+{v_2}="2",
!{(1,-1)}*+{v_3}="3",
!{(1,1)}*+{v_4}="4",
!{(1.5,0)}*+{\scriptstyle{\epsilon_1}},
!{(0.3,-0.3)}*+{\scriptstyle{\epsilon_2}},
"1":"2" "1":"3" "2":"4" "3":"4" "1":"4"
},
$$
then
$$
G/\epsilon_1=
\xygraph{
!{<0pt,0pt>;<10pt,10pt>:<-10pt,10pt>::}
!{(-1,-1)}*+{v_1}="1",
!{(-1,1)}*+{v_2}="2",
!{(1,1)}*+{v_{34}}="3",
"1":"2" "1":"3" "2":"3"
}
\hspace{1cm}\text{and}\hspace{1cm}
G/\epsilon_2=
\xygraph{
!{<0pt,0pt>;<10pt,10pt>:<-10pt,10pt>::}
!{(-1.5,1.5)}*+{v_2}="1"
!{(0,0)}*+{v_{14}}="2"
!{(1.5,-1.5)}*+{v_3}="3"
"2":@/^0.7pc/"1" "1":@/^0.7pc/"2"
"2":@/_0.7pc/"3" "3":@/_0.7pc/"2"
},
$$
so that $\epsilon_1$ is admissible, but $\epsilon_2$ is \emph{not}. Here,
$$
v_{ij}=v_i\sqcup v_j\smallsetminus\{\text{flags constituting edges between $v_i$ and $v_j$}\}.
$$
\eex

\bde
A \emph{contraction sequence} for a graph $G\in\cG_k$ is a sequence of $k-1$ admissible thick edges $(\epsilon_1\in\mathbf{te}(G),\epsilon_2\in\mathbf{te}(G/\epsilon_1),\epsilon_3\in\mathbf{te}((G/\epsilon_1)/\epsilon_2),\ldots)$ The set of contraction sequences for $G$ is denoted $S_G$.
\ede

The set $S_G$ describes different ways of contracting $G$ to a corolla. Every contraction sequence gives rise to a sequence of graphs,
$$
G=G_0\Ra G_1\Ra\cdots\Ra G_{k-2}\Ra G_{k-1},
$$
where $G_i$ has $k-i$ vertices, and is obtained from $G_{i-1}$ by contracting a thick edge. Hence, each $G_i$, for $i\geq 1$, has a unique new vertex $\nu_i$ replacing two old vertices of $G_{i-1}$. To any $s\in S_G$ we associate a levelled binary tree whose legs are labeled by the set of vertices of $G$, and whose $k-1$ vertices are labeled by the vertices $\nu_i$ as in the following example.

\bex
Let $s=(\epsilon_1,\epsilon_2,\epsilon_3)$ be the following contraction sequence
$$
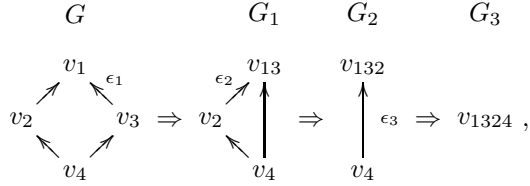
\xygraph{
!{<0pt,0pt>;<20pt,0pt>:<0pt,20pt>::}
!{(0,2)}*+{G},
!{(0,1)}*+{v_1}="1"
!{(-1,0)}*+{v_2}="2"
!{(0.75,0.75)}*+{\scriptstyle{\epsilon_1}}
!{(1,0)}*+{v_3}="3"
!{(0,-1)}*+{v_4}="4"
"4":"2" "4":"3" "2":"1" "3":"1"
}
\Ra
\xygraph{
!{<0pt,0pt>;<20pt,0pt>:<0pt,20pt>::}
!{(0,2)}*+{G_1},
!{(0,1)}*+{v_{13}}="1"
!{(-0.75,0.75)}*+{\scriptstyle{\epsilon_2}}
!{(-1,0)}*+{v_2}="3"
!{(0,-1)}*+{v_4}="4"
"4":"3" "4":"1" "3":"1"
}
\Ra
\xygraph{
!{<0pt,0pt>;<20pt,0pt>:<0pt,20pt>::}
!{(0,2)}*+{G_2},
!{(0,1)}*+{v_{132}}="1"
!{(0.5,0)}*+{\scriptstyle{\epsilon_3}}
!{(0,-1)}*+{v_4}="4"
"4":"1"
}
\Ra
\xygraph{
!{<0pt,0pt>;<20pt,0pt>:<0pt,20pt>::}
!{(0,2)}*+{G_3},
!{(0,0)}*+{v_{1324}}
},
$$
then $\nu_1=v_{13}$, $\nu_2=v_{132}$, $\nu_3=v_{1324}$ and the binary tree associated to $s$ is
$$
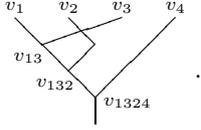
\xygraph{
!{<0pt,0pt>;<10pt,10pt>:<-10pt,10pt>::}
!{(-0.3,2.7)}*+{\scriptstyle{v_1}},
!{(0.7,1.7)}*+{\scriptstyle{v_2}},
!{(1.7,0.7)}*+{\scriptstyle{v_3}},
!{(2.7,-0.3)}*+{\scriptstyle{v_4}},
!{(-1,1.5)}*+{\scriptstyle{v_{13}}},
!{(-1,0.5)}*+{\scriptstyle{v_{132}}},
!{(0,-1.2)}*+{\scriptstyle{v_{1324}}},
!{(-1,-1)}-!{(-0.5,-0.5)},
!{(-0.5,-0.5)}-!{(-0.5,2.5)},
!{(-0.5,0.5)}-!{(0.5,0.5)},
!{(0.5,0.5)}-!{(0.5,1.5)},
!{(-0.5,1.5)}-!{(1.5,0.5)},
!{(-0.5,-0.5)}-!{(2.5,-0.5)},
}.
$$
\eex

On the set $S_G$ we introduce an equivalence relation $\sim$, where $s\sim s'$ if the trees associated to $s$ and $s'$ are isomorphic. This corresponds to forgetting the levels.

\bex
Both of the contraction sequences
$$
s:
\xygraph{
!{<0pt,0pt>;<20pt,0pt>:<0pt,20pt>::}
!{(0,1)}*+{v_1}="1"
!{(-1,0)}*+{v_2}="2"
!{(-0.75,0.75)}*+{\scriptstyle{\epsilon_1}}
!{(1,0)}*+{v_3}="3"
!{(0,-1)}*+{v_4}="4"
"4":"2" "4":"3" "2":"1" "3":"1"
}
\Ra
\xygraph{
!{<0pt,0pt>;<20pt,0pt>:<0pt,20pt>::}
!{(0,1)}*+{v_{12}}="1"
!{(0.75,-0.75)}*+{\scriptstyle{\epsilon_2}}
!{(1,0)}*+{v_3}="3"
!{(0,-1)}*+{v_4}="4"
"4":"3" "4":"1" "3":"1"
}
\Ra
\xygraph{
!{<0pt,0pt>;<20pt,0pt>:<0pt,20pt>::}
!{(0,1)}*+{v_{12}}="1"
!{(0.5,0)}*+{\scriptstyle{\epsilon_3}}
!{(0,-1)}*+{v_{34}}="4"
"4":"1"
}
\Ra
\xygraph{
!{<0pt,0pt>;<20pt,0pt>:<0pt,20pt>::}
!{(0,0)}*+{v_{1234}}
}
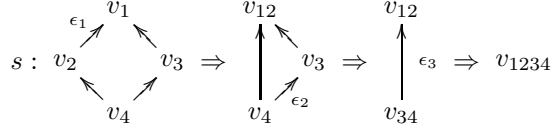
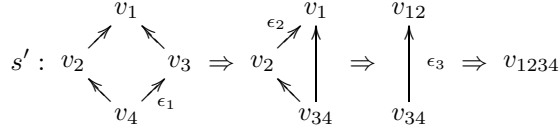
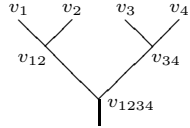
$$
and
$$
s':
\xygraph{
!{<0pt,0pt>;<20pt,0pt>:<0pt,20pt>::}
!{(0,1)}*+{v_1}="1"
!{(-1,0)}*+{v_2}="2"
!{(0.75,-0.75)}*+{\scriptstyle{\epsilon_1}}
!{(1,0)}*+{v_3}="3"
!{(0,-1)}*+{v_4}="4"
"4":"2" "4":"3" "2":"1" "3":"1"
}
\Ra
\xygraph{
!{<0pt,0pt>;<20pt,0pt>:<0pt,20pt>::}
!{(0,1)}*+{v_1}="1"
!{(-0.75,0.75)}*+{\scriptstyle{\epsilon_2}}
!{(-1,0)}*+{v_2}="2"
!{(0,-1)}*+{v_{34}}="3"
"3":"1" "3":"2" "2":"1"
}
\Ra
\xygraph{
!{<0pt,0pt>;<20pt,0pt>:<0pt,20pt>::}
!{(0,1)}*+{v_{12}}="1"
!{(0.5,0)}*+{\scriptstyle{\epsilon_3}}
!{(0,-1)}*+{v_{34}}="4"
"4":"1"
}
\Ra
\xygraph{
!{<0pt,0pt>;<20pt,0pt>:<0pt,20pt>::}
!{(0,0)}*+{v_{1234}}
}
$$
give the tree
$$
\xygraph{
!{<0pt,0pt>;<10pt,10pt>:<-10pt,10pt>::}
!{(-0.3,2.7)}*+{\scriptstyle{v_1}},
!{(0.7,1.7)}*+{\scriptstyle{v_2}},
!{(1.7,0.7)}*+{\scriptstyle{v_3}},
!{(2.7,-0.3)}*+{\scriptstyle{v_4}},
!{(-1,1.5)}*+{\scriptstyle{v_{12}}},
!{(1.5,-1)}*+{\scriptstyle{v_{34}}},
!{(0,-1.2)}*+{\scriptstyle{v_{1234}}},
!{(-1,-1)}-!{(-0.5,-0.5)},
!{(-0.5,-0.5)}-!{(-0.5,2.5)},
!{(1.5,-0.5)}-!{(1.5,0.5)},
!{(-0.5,1.5)}-!{(0.5,1.5)},
!{(-0.5,-0.5)}-!{(2.5,-0.5)},
}
$$
and hence, $s\sim s'$.
\eex

It is not hard to see that if $s$ and $s'$ only differ in the order in which two thick edges in some $G_i$, connecting four distinct vertices, are contracted, then $s\sim s'$. Conversely, if $s\sim s'$, then there is a sequence of contraction sequences $s=s_0,s_1,\ldots,s_{n-1},s_n=s'$ such that, for $0\leq j<n$, $s_j$ and $s_{j+1}$ only differ in the order of which two thick edges in some $G_i$, connecting four distinct vertices, are contracted.

\bde
We denote by $T_G$ the set $S_G/\!\!\sim$ of equivalence classes of contraction sequences.
\ede

Thus, $T_G$ can be identified with a subset of the set of binary trees whose leaves are labeled by the vertices of $G$. From now on we think of elements of $T_G$ as such trees.

If $\epsilon,\epsilon'\in\mathbf{te}(G)$ are two distinct admissible thick edges, then obviously $(G/\epsilon)/\epsilon'=(G/\epsilon')/\epsilon=G/\{\epsilon,\epsilon'\}$. Thus, the notion of contracting a subgraph is well-defined.

\bde
If $H$ is a subgraph of $G$, then $G/H$ denotes the graph where all edges of $H$ have been contracted.
\ede

\bde
We say that a connected subgraph $H$ of $G$ is \emph{admissible} if $G/H$ contains no directed cycles.
\ede

Note that if $H$ is an admissible subgraph of $G$, then every element of $T_H$ describes a way of contracting $G$ to $G/H$.

\subsection{A lemma on trees}\hspace{1pt}

\ble\label{lem:one}
For any graph $G$, given a tree $t\in T_G$ and an internal edge $\epsilon\in\mathbf{e}(t)$, there is a unique tree $t'\in T_G$, different from $t$, and a unique internal edge $\epsilon'\in\mathbf{e}(t')$ such that $t/\epsilon=t'/\epsilon'$.
\ele
\bpro
Consider the vicinity of $\epsilon$ in $t$
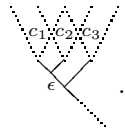
$$
\xygraph{
!{<0pt,0pt>;<5pt,5pt>:<-5pt,5pt>::}
!{(-0.5,0.5)}*+{\scriptstyle{\epsilon}},
!{(1,3)}*+{\scriptstyle{c_1}},
!{(2,2)}*+{\scriptstyle{c_2}},
!{(3,1)}*+{\scriptstyle{c_3}},
!{(0,0)}-!{(2,0)},
!{(0,0)}-!{(0,2)},
!{(0,1)}-!{(1,1)},
!{(0,0)}-!{(0,-1)},
!{(0,-1)}-@{.}!{(0,-3)},
!{(0,2)}-@{.}!{(1,5)},
!{(0,2)}-@{.}!{(3,3)},
!{(1,1)}-@{.}!{(4,2)},
!{(1,1)}-@{.}!{(2,4)},
!{(2,0)}-@{.}!{(3,3)},
!{(2,0)}-@{.}!{(5,1)}
}.
$$

The dashed subtrees describe ways of contracting three admissible disjoint subgraphs of $G$ to corollas. We label them $c_1$, $c_2$ and $c_3$. Two of these are connected by an admissible thick edge and, without loss of generality, we may assume $c_1\la c_2$. Thus, these corollas sit inside $G$ in one of the six following patterns:

$$
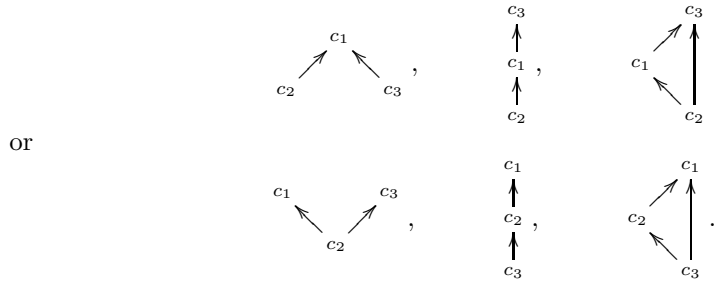
\xygraph{
!{<0pt,0pt>;<20pt,0pt>:<0pt,20pt>::}
!{(-1,-0.5)}*+{\scriptstyle{c_2}}="2"
!{(0,0.5)}*+{\scriptstyle{c_1}}="1"
!{(1,-0.5)}*+{\scriptstyle{c_3}}="3"
"2":"1" "3":"1"
},\hspace{1cm}
\xygraph{
!{<0pt,0pt>;<20pt,0pt>:<0pt,20pt>::}
!{(0,-1)}*+{\scriptstyle{c_2}}="2"
!{(0,0)}*+{\scriptstyle{c_1}}="1"
!{(0,1)}*+{\scriptstyle{c_3}}="3"
"2":"1" "1":"3"
},\hspace{1cm}
\xygraph{
!{<0pt,0pt>;<20pt,0pt>:<0pt,20pt>::}
!{(1,-1)}*+{\scriptstyle{c_2}}="2"
!{(0,0)}*+{\scriptstyle{c_1}}="1"
!{(1,1)}*+{\scriptstyle{c_3}}="3"
"2":"1" "2":"3" "1":"3"
}
$$
or
$$
\xygraph{
!{<0pt,0pt>;<20pt,0pt>:<0pt,20pt>::}
!{(0,-0.5)}*+{\scriptstyle{c_2}}="2"
!{(-1,0.5)}*+{\scriptstyle{c_1}}="1"
!{(1,0.5)}*+{\scriptstyle{c_3}}="3"
"2":"1" "2":"3"
},\hspace{1cm}
\xygraph{
!{<0pt,0pt>;<20pt,0pt>:<0pt,20pt>::}
!{(0,0)}*+{\scriptstyle{c_2}}="2"
!{(0,1)}*+{\scriptstyle{c_1}}="1"
!{(0,-1)}*+{\scriptstyle{c_3}}="3"
"2":"1" "3":"2"
},\hspace{1cm}
\xygraph{
!{<0pt,0pt>;<20pt,0pt>:<0pt,20pt>::}
!{(0,0)}*+{\scriptstyle{c_2}}="2"
!{(1,1)}*+{\scriptstyle{c_1}}="1"
!{(1,-1)}*+{\scriptstyle{c_3}}="3"
"2":"1" "3":"1" "3":"2"
}.
$$

In the first three cases, the required tree has the form

$$
t'=
\xygraph{
!{<0pt,0pt>;<5pt,5pt>:<-5pt,5pt>::}
!{(-0.5,0.5)}*+{\scriptstyle{\epsilon'}},
!{(1,3)}*+{\scriptstyle{c_1}},
!{(2,2)}*+{\scriptstyle{c_2}},
!{(3,1)}*+{\scriptstyle{c_3}},
!{(0,0)}-!{(0,2)},
!{(0,1)}-!{(2,0)},
!{(1,0)}-!{(1,1)},
!{(0,0)}-!{(1,0)},
!{(0,0)}-!{(0,-1)},
!{(0,-1)}-@{.}!{(0,-3)},
!{(0,2)}-@{.}!{(1,5)},
!{(0,2)}-@{.}!{(3,3)},
!{(1,1)}-@{.}!{(4,2)},
!{(1,1)}-@{.}!{(2,4)},
!{(2,0)}-@{.}!{(3,3)},
!{(2,0)}-@{.}!{(5,1)}
}.
$$

In the remaining three cases we have

$$
t'=
\xygraph{
!{<0pt,0pt>;<5pt,5pt>:<-5pt,5pt>::}
!{(1,-0.5)}*+{\scriptstyle{\epsilon'}},
!{(3,1)}*+{\scriptstyle{c_3}},
!{(2,2)}*+{\scriptstyle{c_2}},
!{(1,3)}*+{\scriptstyle{c_1}},
!{(0,0)}-!{(2,0)},
!{(0,0)}-!{(0,2)},
!{(1,0)}-!{(1,1)},
!{(0,0)}-!{(0,-1)},
!{(0,-1)}-@{.}!{(0,-3)},
!{(0,2)}-@{.}!{(1,5)},
!{(0,2)}-@{.}!{(3,3)},
!{(1,1)}-@{.}!{(4,2)},
!{(1,1)}-@{.}!{(2,4)},
!{(2,0)}-@{.}!{(3,3)},
!{(2,0)}-@{.}!{(5,1)}
}.
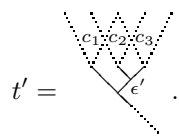
$$
\epro

\subsection{$\Sigma$-bimodules and decorated graphs}\hspace{1pt}

Let $\Sigma$ be the skeleton of the category of finite sets and bijections, with objects $[n]=\{1,\ldots,n\}$. We denote by $\Sigma_n$ the symmetric group of morphisms, $\Sigma([n],[n])$.

\bde
A \emph{($\mathbb{K}$-linear dg) $\Sigma$-bimodule} is a functor from $\Sigma\times\Sigma^{\text{op}}$ to the category of differential graded $\mathbb{K}$-modules.
\ede

That is, a $\Sigma$-bimodule $\cE$ is a collection $\{(\cE(m,n),d_\cE(m,n))\}$ of dg $\mathbb{K}$-vector spaces, such that each $\cE(m,n)$ is endowed with a left $\Sigma_m$-action and a right $\Sigma_n$-action, and these commute. We denote the differential $d_\cE$, dropping the reference to the component. A morphism of $\Sigma$-bimodules $f:\cE\ra\cE'$ is a natural transformation, or more concretely, a collection $\{f{(m,n)}:\cE(m,n)\ra\cE'(m,n)\}$ of degree zero homomorphisms, each commuting with differentials and equivariant with respect to $\Sigma$-actions.

For any finite sets $X$ and $Y$, with $|X|=m$ and $|Y|=n$, we let
$$
\cE(X,Y)=\Bij(X,[m])\times_{\Sigma_m}\cE(m,n)\times_{\Sigma_n}\Bij([n],Y).
$$

\bde
A graph $G$ \emph{decorated} by a $\Sigma$-bimodule $\cE$ is a pair $(G,\otimes_{v\in\mathbf{v}(G)}e_v)$, where
$$
\otimes_{v\in\mathbf{v}(G)}e_v\in\bigotimes_{v\in\mathbf{v}(G)}\cE(out(v),in(v)).
$$
We denote the space of decorations of $G$ by $\cE$ with $\decor{G}{\cE}$.
\ede

Note that the tensor product above is unordered. If $\{V_\alpha\}_{\alpha\in I}$ is a family of $n$ vector spaces, then
$$
\bigotimes_{\alpha\in I}V_\alpha=\left(\bigoplus_{h\in\Bij([n],I)} V_{h(1)}\otimes\ldots\otimes V_{h(n)}\right)_{\Sigma_n},
$$
where $\Sigma_n$ acts by permuting the factors. If $G$ is an $(m,n)$-graph, then $\decor{G}{\cE}$ has natural $\Sigma_m$- and $\Sigma_n$-actions by relabeling.

Let $\cF\cE=\bigoplus_{G\in\cG}\decor{G}{\cE}$, with differential, for $G\in\cG_k$, defined by
\begin{align*}
d_{\cF\cE}&(G,[\otimes_{i=1}^k e_i])=\sum_{i=1}^k(-1)^{|e_1|+\cdots+|e_{i-1}|}(G,[e_1\otimes\cdots\otimes d_\cE e_i\otimes\cdots\otimes e_k]).
\end{align*}
Relabeling gives $\cF\cE$ a structure of $\Sigma$-bimodule, it will be the underlying $\Sigma$-bimodule of the free properad on $\cE$. We note that $\cF\cE$ is weighted, $\cF\cE=\bigoplus_ {k\geq 0}\cF_k\cE$, where $\cF_k\cE=\bigoplus_{G\in\cG_k}\decor{G}{\cE}$. We write $\bar{\cF}\cE=\bigoplus_{k\geq 1}\cF_k\cE$.

\bde
Assume $\cE$ is a $\Sigma$-bimodule. A degree zero morphism of $\Sigma$-bimodules
$$
\mu:\cF_2\cE\ra\cE
$$
we call a \emph{composition}. We write $\mu(e_1,e_2)$ for $\mu(v_1\la v_2,[e_1\otimes e_2])$.
\ede

A composition $\mu$ induces a morphism $\mu_\epsilon:\decor{G}{\cE}\ra\decor{G/\epsilon}{\cE}$ for any admissible thick edge $v_i\xla{\epsilon}v_j$ in $G$, defined by
$$
\mu_\epsilon(G,[e_i\otimes e_j\otimes\otimes_{\nu\neq i,j}e_\nu])=(G/\epsilon,[\mu(e_i,e_j)\otimes\otimes_{\nu\neq i,j}e_\nu]),
$$
where the factor $\mu(e_i,e_j)$ decorates the new vertex of $G/e$.

Any contraction sequence $s$ for an $(m,n)$-graph $G$ hence defines a morphism $\mu_s:\decor{G}{\cE}\ra\cE(m,n)$. The fact that two equivalent sequences induce the same morphism is a consequence of the following lemma.

\ble\label{lem:zero}
Let $v_i,v_j,v_k,v_l$ be four distinct vertices in $G$, and $v_i\xla{\epsilon}v_j$, and $v_k\xla{\epsilon'}v_l$ be contractible thick edges, then $\mu_\epsilon\mu_{\epsilon'}=\mu_{\epsilon'}\mu_\epsilon:\decor{G}{\cE}\ra\decor{G/\{\epsilon,\epsilon'\}}{\cE}$.
\ele
\bpro
We have
\begin{align*}
\mu_\epsilon&\mu_{\epsilon'}(G,[e_k\otimes e_l\otimes e_i\otimes e_j\otimes\otimes_{\nu\neq i,j,k,l}e_\nu]) \\
&=\mu_\epsilon(G/\epsilon',[\mu(e_k,e_l)\otimes e_i\otimes e_j\otimes\otimes_{\nu\neq i,j,k,l}e_\nu]) \\
&=(-1)^{(|e_k|+|e_l|)(|e_i|+|e_j|)}\mu_\epsilon(G/\epsilon',[e_i\otimes e_j\otimes\mu(e_k,e_l)\otimes\otimes_{\nu\neq i,j,k,l}e_\nu]) \\
&=(-1)^{(|e_k|+|e_l|)(|e_i|+|e_j|)}(G/\{\epsilon,\epsilon'\},[\mu(e_i,e_j)\otimes\mu(e_k,e_l)\otimes\otimes_{\nu\neq i,j,k,l}e_\nu]),
\end{align*}
and
\begin{align*}
\mu_{\epsilon'}&\mu_\epsilon(G,[e_k\otimes e_l\otimes e_i\otimes e_j\otimes\otimes_{\nu\neq i,j,k,l}e_\nu]) \\
&=(-1)^{(|e_k|+|e_l|)(|e_i|+|e_j|)}\mu_{\epsilon'}\mu_\epsilon(G,[e_i\otimes e_j\otimes e_k\otimes e_l\otimes\otimes_{\nu\neq i,j,k,l}e_\nu]) \\
&=(-1)^{(|e_k|+|e_l|)(|e_i|+|e_j|)}\mu_{\epsilon'}(G/\epsilon,[\mu(e_i,e_j)\otimes e_k\otimes e_l\otimes\otimes_{\nu\neq i,j,k,l}e_\nu]) \\
&=\mu_{\epsilon'}(G/\epsilon,[e_k\otimes e_l\otimes\mu(e_i,e_j)\otimes\otimes_{\nu\neq i,j,k,l}e_\nu]) \\
&=(G/\{\epsilon,\epsilon'\},[\mu(e_k,e_l)\otimes\mu(e_i,e_j)\otimes\otimes_{\nu\neq i,j,k,l}e_\nu]) \\
&=(-1)^{(|e_k|+|e_l|)(|e_i|+|e_j|)}(G/\{\epsilon,\epsilon'\},[\mu(e_i,e_j)\otimes\mu(e_k,e_l)\otimes\otimes_{\nu\neq i,j,k,l}e_\nu]).
\end{align*}
\epro
Note that the proof relies only on the fact that $\mu$ is of degree zero.

\bde
Let $\cE$ be a $\Sigma$-bimodule, $\mu$ a composition and $t\in T_G$. We denote by $\mu_t$ the degree zero morphism $\decor{G}{\cE}\ra\cE$ induced by a contraction sequence of the equivalence class $t$.
\ede

\bex
If
$$
s:
G=
\xygraph{
!{<0pt,0pt>;<20pt,0pt>:<0pt,20pt>::}
!{(0,1)}*+{v_{1}}="1"
!{(0.75,-0.75)}*+{\scriptstyle{\epsilon_1}}
!{(1,0)}*+{v_2}="3"
!{(0,-1)}*+{v_3}="4"
"4":"3" "4":"1" "3":"1"
}
\Ra
\xygraph{
!{<0pt,0pt>;<20pt,0pt>:<0pt,20pt>::}
!{(0,1)}*+{v_{1}}="1"
!{(0.5,0)}*+{\scriptstyle{\epsilon_2}}
!{(0,-1)}*+{v_{23}}="4"
"4":"1"
}
\Ra
\xygraph{
!{<0pt,0pt>;<20pt,0pt>:<0pt,20pt>::}
!{(0,0)}*+{v_{123}}
},
$$
then
$$
t=
\xygraph{
!{<0pt,0pt>;<5pt,5pt>:<-5pt,5pt>::}
!{(-0.2,2.8)}*+{\scriptstyle{v_1}},
!{(1.3,1.3)}*+{\scriptstyle{v_2}},
!{(2.8,-0.2)}*+{\scriptstyle{v_3}},
!{(-1,-1)}-!{(-0.5,-0.5)},
!{(-0.5,-0.5)}-!{(-0.5,2.5)},
!{(1,-0.5)}-!{(1,1)},
!{(-0.5,-0.5)}-!{(2.5,-0.5)},
}
$$
and for $(G,[e_1\otimes e_2\otimes e_3])\in\decor{G}{\cE}$ we have
\begin{align*}
\mu_t(G,[e_1\otimes e_2\otimes e_3])=&\mu_{\epsilon_2}\mu_{\epsilon_1}(G,[e_1\otimes e_2\otimes e_3]) \\
=&(-1)^{|e_1|(|e_2|+|e_3|)}\mu_{\epsilon_2}\mu_{\epsilon_1}(G,[e_2\otimes e_3\otimes e_1]) \\
=&(-1)^{|e_1|(|e_2|+|e_3|)}\mu_{\epsilon_2}(G/\epsilon_1,[\mu(e_2,e_3)\otimes e_1]) \\
=&\mu_{\epsilon_2}(G/\epsilon_1,[e_1\otimes\mu(e_2,e_3)])=\mu(e_1,\mu(e_2,e_3)).
\end{align*}
\eex

\subsection{Properads and bar construction}\hspace{1pt}

A composition $\mu$, such that for any graph $G$ and $t,t'\in T_G$, the two morphisms $\mu_t$ and $\mu_{t'}$ are equal, is called \emph{associative}. A composition is associative if and only if it is associative on graphs with three vertices.

Given a $\Sigma$-bimodule $\cP$ and a composition $\mu$, a morphism $d:\cP\ra\cP$ is a \emph{derivation} of $\cP$ with respect to $\mu$ if $d\mu=\mu(d,\Id)+\mu(\Id,d):\cF_2\cP\ra\cP$.

\bde
A (dg) \emph{properad without unit} is a pair $(\cP,\mu)$, where $\cP$ is a dg $\Sigma$-bimodule and $\mu$ is an associative composition such that $d_\cP$ is a derivation with respect to $\mu$. If $\eta:\mathbb{K}\ra\cP(1,1)$ is such that $\mu(\eta(1),p)=\mu(p,\eta(1))=p$ for all $p$ we say that $\eta$ is a \emph{unit} for $\mu$. We call a triple $(\cP,\mu,\eta)$ a \emph{properad}.
\ede

\bex
With $\mu$ as grafting legs of graphs $\bar{\cF}\cE$ is the \emph{free properad without unit} on the $\Sigma$-bimodule $\cE$. With the trivial graph as unit, $\cF\cE$ is the \emph{free properad with unit} on $\cE$.
\eex

Any derivation $d$ of $\bar{\cF}\cE$ is determined by its restriction, $d|_\cE:\cE\ra\bar{\cF}\cE$, to generators. If the image of $d|_\cE$ is in $\oplus_{i\leq k}\bar{\cF}_i\cE$ we say that $d$ is of order $k$. Hence, $d_\cE$ induces the linear, or order $1$, derivation $d_{\cF\cE}$. For $k=2$ we say that $d$ is quadratic.

\bde
Given a $\Sigma$-bimodule $\cE=(\{\cE(m,n)\},d_\cE)$, where $\cE(m,n)=\oplus_i\cE(m,n)^i$, the \emph{shifted by $j$} $\Sigma$-bimodule $\cE[j]$ is defined by $\cE[j](m,n)^i=\cE(m,n)^{i+j}$.
\ede

Consider $\cE[1]$. We let $s^{-1}$ be the degree $-1$ map identifying $\cE^i$ with $\cE[1]^{i-1}$, and write elements of $\cE[1]^{i-1}$ as $s^{-1}e$ with $e\in\cE^i$. The differential $d_\cE$ on $\cE$ induces a differential on $\cE[1]$ by $d_{\cE[1]}s^{-1}e=-s^{-1}d_\cE e$. We usually identify $d_{\cE[1]}$ with $d_{\cE}$, remembering the anti-commutation law with $s^{-1}$. In the same manner we write $d_{\cF\cE}$ for $d_{\cF(\cE[1])}$.

\bde
A \emph{coproperad} is a $\Sigma$-bimodule $\cC$ together with a quadratic differential $\delta$ on $\bar{\cF}(\cC[-1])$, that is, a degree $1$ derivation $\delta:\bar{\cF}(\cC[-1])\ra\bar{\cF}(\cC[-1])$ such that the image of $\delta|_{\cC[-1]}$ is in $\cC[-1]\oplus\bar{\cF}_2(\cC[-1])$ and $\delta^2=0$.
\ede

Given a coproperad $(\cC,\delta)$, the morphism
$$
\tilde{\Delta}:\cC[-1]\xra{\delta|_{\cC[-1]}}\bar{\cF}\cC[-1]\xra{proj}\cF_2\cC[-1]
$$
induces a degree zero morphism $\Delta:\cC\ra\bar{\cF}_2\cC$ called the \emph{cocomposition}.

Consider the degree $1$ morphism $\tilde{\Delta}:(\bar{\cF}\cE)[-1]\ra\bar{\cF}_2((\bar{\cF}\cE)[-1])$,
\begin{align*}
\tilde{\Delta}(G,&s\otimes_{v\in\mathbf{v}(G)}e_v) \\
=&\sum_{G=G_1\circ G_2}(-1)^\varepsilon(v_1\la v_2,[(G_1,s\otimes_{v\in\mathbf{v}(G_1)}e_v)\otimes(G_2,s\otimes_{v\in\mathbf{v}(G_2)}e_v)])
\end{align*}
where $\varepsilon$ is induced by the Koszul convention and also includes a term $\sum_{v\in\mathbf{v}(G_1)}|e_v|$. It determines a derivation $\delta_\Delta$ on $\bar{\cF}((\bar{\cF}\cE)[-1])$, and it is a fact that the quadratic degree $1$ derivation $\delta=d_{\cF\cE}+\delta_\Delta$ is a differential on $\bar{\cF}((\bar{\cF}\cE)[-1])$. This gives $\bar{\cF}\cE$ a structure of coproperad, it is called the \emph{cofree} coproperad on $\cE$, and equipped with this structure we denote it $\bar{\cF}^c\cE$. The degree zero morphism $\Delta:\bar{\cF}^c\cE\ra\bar{\cF}_2\bar{\cF}^c\cE$ is not coassociative in the sense that $(\Delta,\Id)\Delta=(\Id,\Delta)\Delta$ as the following example shows.

\bex
Let
$$
G=
\xygraph{
!{<0pt,0pt>;<10pt,0pt>:<0pt,10pt>::}
!{(-2,-1)}*+{\scriptstyle{v_2}}="2"
!{(0,1)}*+{\scriptstyle{v_1}}="1"
!{(2,-1)}*+{\scriptstyle{v_3}}="3"
"2":"1" "3":"1"
}.
$$
Then
\begin{align*}
(\Delta,\Id)&\Delta(G,[e_1\otimes e_2\otimes e_3]) \\
=&(\Delta,\Id)(((v_1\la v_2)\la v_3,[(v_1\la v_2,[e_1\otimes e_2])\otimes(v_3,[e_3])]) \\
&+(-1)^{|e_2||e_3|}((v_1\la v_3)\la v_2,[(v_1\la v_3,[e_1\otimes e_3])\otimes(v_2,[e_2])])) \\
=&(G,[(v_1,[e_1])\otimes(v_2,[e_2])\otimes(v_3,[e_3])]) \\
&+(-1)^{|e_2||e_3|}(G,[(v_1,[e_1])\otimes(v_3,[e_3])\otimes(v_2,[e_2])]) \\
=&2(G,[(v_1,[e_1])\otimes(v_2,[e_2])\otimes(v_3,[e_3])]),
\end{align*}
while $(\Id,\Delta)\Delta(G,[e_1\otimes e_2\otimes e_3])=0$. However, $\delta$ is still a differential as indicated by
\begin{align*}
(\tilde{\Delta},\Id)&\tilde{\Delta}(G,s[e_1\otimes e_2\otimes e_3]) \\
=&(\tilde{\Delta},\Id)((-1)^{|e_1|+|e_2|}((v_1\la v_2)\la v_3,[(v_1\la v_2,s[e_1\otimes e_2])\otimes(v_3,s[e_3])]) \\
&+(-1)^{|e_2||e_3|+|e_1|+|e_3|}((v_1\la v_3)\la v_2,[(v_1\la v_3,s[e_1\otimes e_3])\otimes(v_2,s[e_2])]) \\
=&(-1)^{|e_2|}(G,[(v_1,s[e_1])\otimes(v_2,s[e_2])\otimes(v_3,s[e_3])]) \\
&+(-1)^{|e_2||e_3|+|e_3|}(G,[(v_1,s[e_1])\otimes(v_3,s[e_3])\otimes(v_2,s[e_2])]) \\
=&((-1)^{|e_2|}+(-1)^{|e_2|+1})(G,[(v_1,s[e_1])\otimes(v_2,s[e_2])\otimes(v_3,s[e_3])])=0.
\end{align*}
\eex

A morphism $\partial:\cC\ra\cC$ is a \emph{coderivation} if $\Delta\partial=(\partial,\Id)\Delta+(\Id,\partial)\Delta$. Any coderivation $\partial$ on $\bar{\cF}^c\cE$ is determined by its projection
$$
\partial_1:\bar{\cF}^c\cE\xra{\partial}\bar{\cF}^c\cE\xra{proj}\cE,
$$
and if $\partial_1$ is zero on $\oplus_{i>k}\bar{\cF}^c_i\cE$ we say that $\partial$ is of order $k$. For $k=2$ we say that $\partial$ is quadratic. If $\cP$ is a properad, the composition $\mu$ induces a degree $1$ morphism $\tilde{\mu}:\bar{\cF}_2(\cP[1])\ra\cP[1]$ by
$$
\tilde{\mu}(s^{-1}p_1,s^{-1}p_2)=(-1)^{|p_1|}s^{-1}\mu(p_1,p_2),
$$
which determines a degree $1$ quadratic coderivation $\partial_\mu$ on $\bar{\cF}^c(\cP[1])$. The fact that $(\cP,\mu)$ is a properad is then equivalent to $d_{\cF\cP}+\partial_\mu$ being a codifferential.

\bde
The \emph{bar construction}, $B\cP$, on a properad $\cP$, is the coproperad \\
$(\bar{\cF}^c(\cP[1]),\partial_\cP)$, with codifferential $\partial_\cP=d_{\cF\cP}+\partial_\mu$.
\ede

Reference for properads and their bar contructions is \cite{Vallette2004}.

\section{Sh properads and transfer}

\subsection{Strong homotopy properads}\hspace{1pt}

\bde
A \emph{strong homotopy (or sh) properad} is a pair $(\cE,\partial_\cE)$, where $\cE$ is a $\Sigma$-bimodule and $\partial_\cE$ is a codifferential on $\bar{\cF}^c(\cE[1])$.
\ede
A morphism of sh properads $F:\cE\ra\cE'$ is a morphism of dg coproperads \\
$\bar{\cF}^c(\cE[1])\ra\bar{\cF}^c(\cE'[1])$.

\subsection{Transfer}\hspace{1pt}

We now turn our attention to the following situation. Let $\cP$ be a dg properad, and $\cE$ a dg $\Sigma$-bimodule. Assume we have degree zero morphisms $g:\cP\ra\cE$ and $f:\cE\ra\cP$ of dg $\Sigma$-bimodules, such that $fg$ is homotopic to the identity on $\cP$, that is, there is a degree $-1$ morphism $h:\cP\ra\cP$ such that
$$
fg-\Id_\cP=d_\cP h+h d_\cP.
$$
Out of these data we want to construct a codifferential $\partial_\cE$ on $\bar{\cF}^c(\cE[1])$, and hence get explicit formulae for an induced structure of sh properad on $\cE$.

The morphisms $f,g$ and $h$ induce morphisms $\cE[1]\ra\cP[1]$, $\cP[1]\ra\cE[1]$ and $\cP[1]\ra\cP[1]$ denoted by the same letters. Note that $hs^{-1}p=-s^{-1}hp$ for any $p\in\cP$.

Composing the morphism $\tilde{\mu}$ with $h$ we get a degree zero morphism $\cF_2(\cP[1])\ra\cP[1]$, pictorially represented as
$$
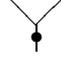
\xygraph{
!{<0pt,0pt>;<10pt,10pt>:<-10pt,10pt>::}
!{(0,0)}-!{(1,0)},
!{(0,0)}-!{(0,1)},
!{(0,0)}-!{(-0.5,-0.5)},
!{(-0.25,-0.25)}*{\bullet}
}
$$
with $\tilde{\mu}$ decorating the vertex and $h$ decorating the dot (see \cite{Kontsevich2001}). As the proof of Lemma \ref{lem:zero} shows, for any graph $G$ with $|\mathbf{v}(G)|\geq 2$ and any $t\in T_G$, there is to $h\tilde{\mu}$ an associated degree zero morphism $(h\tilde{\mu})_t:\decor{G}{\cP[1]}\ra\cP[1]$. This in turn determines a degree $1$ morphism $\theta_t:\decor{G}{\cP[1]}\ra\cP[1]$, such that $(h\tilde{\mu})_t=h\theta_t$.

Composing $\theta_t$ with $fg$ we get a second morphism $\theta_t^\circ:\decor{G}{\cP[1]}\ra\cP[1]$. We represent $fg$ pictorially as $\circ$.

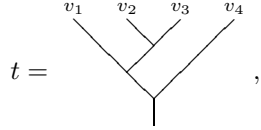
\bex
If
$$
t=
\xygraph{
!{<0pt,0pt>;<10pt,10pt>:<-10pt,10pt>::}
!{(-0.3,2.7)}*+{\scriptstyle{v_1}},
!{(0.7,1.7)}*+{\scriptstyle{v_2}},
!{(1.7,0.7)}*+{\scriptstyle{v_3}},
!{(2.7,-0.3)}*+{\scriptstyle{v_4}},
!{(-1,-1)}-!{(-0.5,-0.5)},
!{(-0.5,-0.5)}-!{(-0.5,2.5)},
!{(-0.5,0.5)}-!{(1.5,0.5)},
!{(0.5,0.5)}-!{(0.5,1.5)},
!{(-0.5,-0.5)}-!{(2.5,-0.5)},
},
$$
then $(h\tilde{\mu})_t$, $\theta_t$ and $\theta_t^\circ$ are represented respectively by
$$
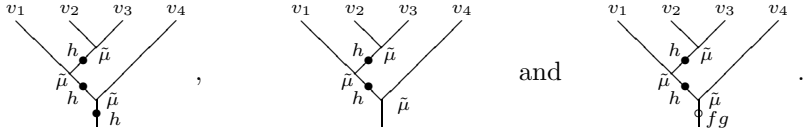
\xygraph{
!{<0pt,0pt>;<10pt,10pt>:<-10pt,10pt>::}
!{(-0.3,2.7)}*+{\scriptstyle{v_1}},
!{(0.7,1.7)}*+{\scriptstyle{v_2}},
!{(1.7,0.7)}*+{\scriptstyle{v_3}},
!{(2.7,-0.3)}*+{\scriptstyle{v_4}},
!{(-0.8,0.5)}*+{\scriptstyle{\tilde{\mu}}},
!{(0.5,0.2)}*+{\scriptstyle{\tilde{\mu}}},
!{(-0.2,-0.8)}*+{\scriptstyle{\tilde{\mu}}},
!{(-0.5,-1.2)}*+{\scriptstyle{h}},
!{(0,0.9)}*+{\scriptstyle{h}},
!{(-0.9,0)}*+{\scriptstyle{h}},
!{(-1,-1)}-!{(-0.5,-0.5)},
!{(-0.5,-0.5)}-!{(-0.5,2.5)},
!{(-0.5,0.5)}-!{(1.5,0.5)},
!{(0.5,0.5)}-!{(0.5,1.5)},
!{(-0.5,-0.5)}-!{(2.5,-0.5)},
!{(0,0.5)}*+{\scriptstyle{\bullet}},
!{(-0.5,0)}*+{\scriptstyle{\bullet}},
!{(-0.75,-0.75)}*+{\scriptstyle{\bullet}},
},
\hspace{1cm}
\xygraph{
!{<0pt,0pt>;<10pt,10pt>:<-10pt,10pt>::}
!{(-0.3,2.7)}*+{\scriptstyle{v_1}},
!{(0.7,1.7)}*+{\scriptstyle{v_2}},
!{(1.7,0.7)}*+{\scriptstyle{v_3}},
!{(2.7,-0.3)}*+{\scriptstyle{v_4}},
!{(-0.8,0.5)}*+{\scriptstyle{\tilde{\mu}}},
!{(0.5,0.2)}*+{\scriptstyle{\tilde{\mu}}},
!{(-0.2,-1)}*+{\scriptstyle{\tilde{\mu}}},
!{(0,0.9)}*+{\scriptstyle{h}},
!{(-0.9,0)}*+{\scriptstyle{h}},
!{(-1,-1)}-!{(-0.5,-0.5)},
!{(-0.5,-0.5)}-!{(-0.5,2.5)},
!{(-0.5,0.5)}-!{(1.5,0.5)},
!{(0.5,0.5)}-!{(0.5,1.5)},
!{(-0.5,-0.5)}-!{(2.5,-0.5)},
!{(0,0.5)}*+{\scriptstyle{\bullet}},
!{(-0.5,0)}*+{\scriptstyle{\bullet}},
}
\hspace{0.5cm}\text{and}\hspace{0.5cm}
\xygraph{
!{<0pt,0pt>;<10pt,10pt>:<-10pt,10pt>::}
!{(-0.3,2.7)}*+{\scriptstyle{v_1}},
!{(0.7,1.7)}*+{\scriptstyle{v_2}},
!{(1.7,0.7)}*+{\scriptstyle{v_3}},
!{(2.7,-0.3)}*+{\scriptstyle{v_4}},
!{(-0.8,0.5)}*+{\scriptstyle{\tilde{\mu}}},
!{(0.5,0.2)}*+{\scriptstyle{\tilde{\mu}}},
!{(-0.2,-0.8)}*+{\scriptstyle{\tilde{\mu}}},
!{(-0.5,-1.2)}*+{\scriptstyle{fg}},
!{(0,0.9)}*+{\scriptstyle{h}},
!{(-0.9,0)}*+{\scriptstyle{h}},
!{(-1,-1)}-!{(-0.5,-0.5)},
!{(-0.5,-0.5)}-!{(-0.5,2.5)},
!{(-0.5,0.5)}-!{(1.5,0.5)},
!{(0.5,0.5)}-!{(0.5,1.5)},
!{(-0.5,-0.5)}-!{(2.5,-0.5)},
!{(0,0.5)}*+{\scriptstyle{\bullet}},
!{(-0.5,0)}*+{\scriptstyle{\bullet}},
!{(-0.75,-0.75)}*+{\scriptstyle{\circ}},
}.
$$
\eex

Given a graph $G$ and $t\in T_G$, any internal edge $\epsilon$ in $t$ determines two subtrees of $t$, one with $\epsilon$ as root and one with $\epsilon$ as leaf. We denote them $t_r$ and $t_l$ respectively. The subtree $t_r$ determines an admissible subgraph $H_\epsilon$ of $G$, consisting of the vertices labeling leaves of $t_r$. Obviously, $t_r\in T_{H_\epsilon}$, implying that $\theta_{t_r}$ and $\theta_{t_r}^\circ$ may be viewed as a morphisms $\decor{G}{\cP[1]}\ra\decor{G/H_\epsilon}{\cP[1]}$. Moreover, $t_l\in T_{G/H_\epsilon}$. Hence, $\epsilon$ determines the following two morphisms $\decor{G}{\cP[1]}\ra\cP[1]$:
\begin{align*}
\theta_{t,\epsilon}^{\Id}:&\decor{G}{\cP[1]}\xra{\theta_{t_r}}\decor{G/H_\epsilon}{\cP[1]}\xra{\theta_{t_l}}\cP[1], \\
\theta_{t,\epsilon}^\circ:&\decor{G}{\cP[1]}\xra{\theta_{t_r}^\circ}\decor{G/H_\epsilon}{\cP[1]}\xra{\theta_{t_l}}\cP[1].
\end{align*}

\bex
If
$$
t=
\xygraph{
!{<0pt,0pt>;<10pt,10pt>:<-10pt,10pt>::}
!{(-0.3,2.7)}*+{\scriptstyle{v_1}},
!{(0.7,1.7)}*+{\scriptstyle{v_2}},
!{(1.7,0.7)}*+{\scriptstyle{v_3}},
!{(2.7,-0.3)}*+{\scriptstyle{v_4}},
!{(-0.8,0)}*+{\scriptstyle{\epsilon}},
!{(-1,-1)}-!{(-0.5,-0.5)},
!{(-0.5,-0.5)}-!{(-0.5,2.5)},
!{(-0.5,0.5)}-!{(1.5,0.5)},
!{(0.5,0.5)}-!{(0.5,1.5)},
!{(-0.5,-0.5)}-!{(2.5,-0.5)},
},
$$
then $\theta_{t,\epsilon}^{\Id}$ and $\theta_{t,\epsilon}^\circ$ are represented respectively by
$$
\xygraph{
!{<0pt,0pt>;<10pt,10pt>:<-10pt,10pt>::}
!{(-0.3,2.7)}*+{\scriptstyle{v_1}},
!{(0.7,1.7)}*+{\scriptstyle{v_2}},
!{(1.7,0.7)}*+{\scriptstyle{v_3}},
!{(2.7,-0.3)}*+{\scriptstyle{v_4}},
!{(-0.8,0.5)}*+{\scriptstyle{\tilde{\mu}}},
!{(0.5,0.2)}*+{\scriptstyle{\tilde{\mu}}},
!{(-0.2,-1)}*+{\scriptstyle{\tilde{\mu}}},
!{(0,0.9)}*+{\scriptstyle{h}},
!{(-1,-1)}-!{(-0.5,-0.5)},
!{(-0.5,-0.5)}-!{(-0.5,2.5)},
!{(-0.5,0.5)}-!{(1.5,0.5)},
!{(0.5,0.5)}-!{(0.5,1.5)},
!{(-0.5,-0.5)}-!{(2.5,-0.5)},
!{(0,0.5)}*+{\scriptstyle{\bullet}},
}
\hspace{2cm}\text{and}\hspace{2cm}
\xygraph{
!{<0pt,0pt>;<10pt,10pt>:<-10pt,10pt>::}
!{(-0.3,2.7)}*+{\scriptstyle{v_1}},
!{(0.7,1.7)}*+{\scriptstyle{v_2}},
!{(1.7,0.7)}*+{\scriptstyle{v_3}},
!{(2.7,-0.3)}*+{\scriptstyle{v_4}},
!{(-0.8,0.5)}*+{\scriptstyle{\tilde{\mu}}},
!{(0.5,0.2)}*+{\scriptstyle{\tilde{\mu}}},
!{(-0.2,-1)}*+{\scriptstyle{\tilde{\mu}}},
!{(0,0.9)}*+{\scriptstyle{h}},
!{(-0.9,0)}*+{\scriptstyle{fg}},
!{(-1,-1)}-!{(-0.5,-0.5)},
!{(-0.5,-0.5)}-!{(-0.5,2.5)},
!{(-0.5,0.5)}-!{(1.5,0.5)},
!{(0.5,0.5)}-!{(0.5,1.5)},
!{(-0.5,-0.5)}-!{(2.5,-0.5)},
!{(0,0.5)}*+{\scriptstyle{\bullet}},
!{(-0.5,0)}*+{\scriptstyle{\circ}},
}.
$$
\eex

We now prove a lemma about the morphisms $\theta_{t,\epsilon}^{\Id}$ needed later.
\ble\label{lem:three}
We have
$$
\sum_{t\in T_G}\sum_{\epsilon\in\mathbf{e}(t)}\theta_{t,\epsilon}^{\Id}=0
$$
as a morphism $\decor{G}{\cP[1]}\ra\cP[1]$.
\ele
\bpro
For any pair $(t,\epsilon)$, where $t\in T_G$ and $\epsilon$ is an internal edge of $t$, there is by Lemma \ref{lem:one} the unique pair $(t',\epsilon')$ such that $t/\epsilon=t'/\epsilon'$. Hence, it is enough to show that $\theta_{t,\epsilon}^{\Id}+\theta_{t',\epsilon'}^{\Id}=0$, which in turn is enough to check for the three-vertex graphs listed in the proof of Lemma \ref{lem:one}. Consider the first graph
$$
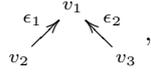
\xygraph{
!{<0pt,0pt>;<20pt,0pt>:<0pt,20pt>::}
!{(-1,-0.5)}*+{\scriptstyle{v_2}}="2"
!{(0,0.5)}*+{\scriptstyle{v_1}}="1"
!{(1,-0.5)}*+{\scriptstyle{v_3}}="3"
!{(-0.75,0.25)}*+{\scriptstyle{\epsilon_1}}
!{(0.75,0.25)}*+{\scriptstyle{\epsilon_2}}
"2":"1" "3":"1"
},
$$
which gives the two trees
$$
t=
\xygraph{
!{<0pt,0pt>;<10pt,10pt>:<-10pt,10pt>::},
!{(0.2,2.3)}*+{\scriptstyle{v_1}},
!{(1.2,1.2)}*+{\scriptstyle{v_2}},
!{(2.3,0.2)}*+{\scriptstyle{v_3}},
!{(-0.3,0.5)}*+{\scriptstyle{\epsilon}},
!{(-0.5,-0.5)}-!{(0,0)},
!{(0,0)}-!{(2,0)},
!{(0,0)}-!{(0,2)},
!{(0,1)}-!{(1,1)}
}
\hspace{1cm}\text{and}\hspace{1cm}
t'=
\xygraph{
!{<0pt,0pt>;<10pt,10pt>:<-10pt,10pt>::},
!{(0.2,2.3)}*+{\scriptstyle{v_1}},
!{(1.2,1.2)}*+{\scriptstyle{v_2}},
!{(2.3,0.2)}*+{\scriptstyle{v_3}},
!{(-0.3,0.5)}*+{\scriptstyle{\epsilon'}},
!{(-0.5,-0.5)}-!{(0,0)},
!{(0,0)}-!{(0,2)},
!{(1,0)}-!{(1,1)},
!{(0,0)}-!{(1,0)},
!{(0,1)}-!{(2,0)},
},
$$
corresponding to
$$
\theta_{t,\epsilon}^{\Id}=\tilde{\mu}_{\epsilon_2}\tilde{\mu}_{\epsilon_1}\hspace{1cm}\text{and}\hspace{1cm}\theta_{t',\epsilon'}^{\Id}=\tilde{\mu}_{\epsilon_1}\tilde{\mu}_{\epsilon_2}.
$$
Now
\begin{align*}
\tilde{\mu}_{\epsilon_2}\tilde{\mu}_{\epsilon_1}(G,[s^{-1}p_1\otimes s^{-1}p_2\otimes s^{-1}p_3])=&\tilde{\mu}_{\epsilon_2}(G/\epsilon_1,[\tilde{\mu}(s^{-1}p_1,s^{-1}p_2)\otimes s^{-1}p_3]) \\
=&\tilde{\mu}(\tilde{\mu}(s^{-1}p_1,s^{-1}p_2),s^{-1}p_3),
\end{align*}
while
\begin{align*}
\tilde{\mu}_{\epsilon_1}\tilde{\mu}_{\epsilon_2}&(G,[s^{-1}p_1\otimes s^{-1}p_2\otimes s^{-1}p_3]) \\
=&(-1)^{|p_2||p_3|+|p_2|+|p_3|+1}\tilde{\mu}_{\epsilon_1}\tilde{\mu}_{\epsilon_2}(G,[s^{-1}p_1\otimes s^{-1}p_3\otimes s^{-1}p_2]) \\
=&(-1)^{|p_2||p_3|+|p_2|+|p_3|+1}\tilde{\mu}_{\epsilon_1}(G/\epsilon_2,[\tilde{\mu}(s^{-1}p_1,s^{-1}p_3)\otimes s^{-1}p_2]) \\
=&(-1)^{|p_2||p_3|+|p_2|+|p_3|+1}\tilde{\mu}(\tilde{\mu}(s^{-1}p_1,s^{-1}p_3),s^{-1}p_2),
\end{align*}
and since $\mu$ is associative
\begin{align*}
\tilde{\mu}(\tilde{\mu}(s^{-1}p_1,s^{-1}p_2),s^{-1}p_3)=&(-1)^{|p_1|}\tilde{\mu}(s^{-1}\mu(p_1,p_2),s^{-1}p_3) \\
=&(-1)^{|p_2|}s^{-1}\mu(\mu(p_1,p_2),p_3) \\
=&(-1)^{|p_2||p_3|+|p_2|}s^{-1}\mu(\mu(p_1,p_3),p_2) \\
=&(-1)^{|p_2||p_3|+|p_1|+|p_2|+|p_3|}\tilde{\mu}(s^{-1}\mu(p_1,p_3),s^{-1}p_2) \\
=&(-1)^{|p_2||p_3|+|p_2|+|p_3|}\tilde{\mu}(\tilde{\mu}(s^{-1}p_1,s^{-1}p_3),s^{-1}p_2),
\end{align*}
showing that $\tilde{\mu}_{\epsilon_2}\tilde{\mu}_{\epsilon_1}+\tilde{\mu}_{\epsilon_1}\tilde{\mu}_{\epsilon_2}=0$. The calculations for the other five cases are similar, so we omit them.
\epro

\bde
For a graph $G$ with $|\mathbf{v}(G)|\geq 2$, we set
$$
\theta_G=\sum_{t\in T_G}\theta_t:\decor{G}{\cP[1]}\ra\cP[1].
$$
Next, we define $\partial_G:\decor{G}{\cE[1]}\ra\cE[1]$ to be $d_\cE$ on corollas, and for $G$ such that $|\mathbf{v}(G)|\geq 2$,
$$
\partial_G=g\theta_G f^{\otimes k}.
$$
\ede

The family $\partial_G$ determines a degree $1$ coderivation $\partial_\cE$ on $\bar{\cF}^c(\cE[1])$. Given an admissible subgraph $H\subset G$, we may view $\partial_H$ as a morphism $\decor{G}{\cE[1]}\ra\decor{G/H}{\cE[1]}$ by decorating the new vertex of $G/H$ by $\partial_H$ applied to the decorated subgraph $H$. Then, on a decorated graph $G$,
$$
\partial_\cE=\sum_{H\subset G}\partial_H,
$$
where $H$ runs over all admissible subgraphs of $G$.

Given a morphism $F:\bar{\cF}^c(\cE)\ra\bar{\cF}^c(\cE')$, we will write $F_k$ for the composition
$$
\bar{\cF}^c(\cE)\xra{F}\bar{\cF}^c(\cE')\xra{proj}\bar{\cF}_k^c(\cE').
$$

\bth\label{thm:one}
Let $\cP$ be a properad and $\cE$ a $\Sigma$-bimodule, $g:\cP\ra\cE$, $f:\cE\ra\cP$ be morphisms of $\Sigma$-bimodules such that $fg-\Id_\cP=d_\cP h+h d_\cP$, where $h:\cP\ra\cP$ is a morphism of degree $-1$. Then the morphism $\partial_\cE$ constructed above is a codifferential on $\bar{\cF}(\cE[1])$.
\eth
\bpro
The equality
\begin{align*}
\Delta\partial_\cE^2=&(\partial_\cE,\Id)\Delta\partial_\cE+(\Id,\partial_\cE)\Delta\partial_\cE \\
=&(\partial_\cE^2,\Id)\Delta+(\partial_\cE,\partial_\cE)\Delta-(\partial_\cE,\partial_\cE)\Delta+(\Id,\partial_\cE^2)\Delta \\
=&(\partial_\cE^2,\Id)\Delta+(\Id,\partial_\cE^2)\Delta,
\end{align*}
shows that $\partial_\cE^2$ is a coderivation, and hence that it is determined by $(\partial_\cE^2)_1$. Thus, we need to show that on any decorated graph $G$ we have
$$
(\partial_\cE^2)_1=\sum_{H\subset G}\partial_{G/H}\partial_H=0.
$$
Rewriting the terms with $H=G$, and $H$ a corolla, according to the definition of $\partial_H$, we see that this is equivalent to
$$
d_\cE g\theta_G f^{\otimes k}+g\theta_G f^{\otimes k}d_{\cF\cE}=-\sum_{H\subset G}\partial_{G/H}\partial_H
$$
where the sum is over proper admissible subgraphs $H$ with at least two vertices. Rewriting the right hand side we see that this is
$$
d_\cE g\theta_G f^{\otimes k}+g\theta_G f^{\otimes k}d_{\cF\cE}=\sum_{t\in T_G}\sum_{\epsilon\in\mathbf{e}(t)}-g\theta_{t,\epsilon}^\circ f^{\otimes k},
$$
which would follow from
\begin{equation}\label{eqn:one}
d_\cP \theta_G+\theta_G d_{\cF\cP}=\sum_{t\in T_G}\sum_{\epsilon\in\mathbf{e}(t)}-\theta_{t,\epsilon}^\circ.
\end{equation}
By Lemma \ref{lem:three}, this is equivalent to
$$
d_\cP \theta_G+\theta_G d_{\cF\cP}=\sum_{t\in T_G}\sum_{\epsilon\in\mathbf{e}(t)}\theta_{t,\epsilon}^{\Id}-\theta_{t,\epsilon}^\circ,
$$
which in turn is equivalent to
$$
\sum_{t\in T_G}\left(d_\cP \theta_t+\theta_t d_{\cF\cP}-\sum_{\epsilon\in\mathbf{e}(t)}(\theta_{t,\epsilon}^{\Id}-\theta_{t,\epsilon}^\circ)\right)=0.
$$
We now show, by induction over $|\mathbf{v}(G)|$, that each summand corresponding to a tree is zero. For $|\mathbf{v}(G)|=2$ it reads
$$
d_\cP \tilde{\mu}+\tilde{\mu}(d_\cP,\Id_\cP)+\tilde{\mu}(\Id_\cP,d_\cP)=0,
$$
which precisely states the fact that $d_\cP$ is a derivation on $\cP[1]$. Now, assume it true for $|\mathbf{v}(G)|<k$. For a fixed tree $t\in T_G$ with non-trivial principal subtrees $t_1, t_2$ we get
\begin{align}\label{eqn:two}\begin{split}
d_\cP \theta_t=&d_\cP \tilde{\mu}(h\theta_{t_1},h\theta_{t_2}) \\
=&-\tilde{\mu}(d_\cP h\theta_{t_1},h\theta_{t_2})-\tilde{\mu}(h\theta_{t_1},d_\cP h\theta_{t_2}) \\
=&-\tilde{\mu}(fg \theta_{t_1},h\theta_{t_2})+\tilde{\mu}(hd_\cP\theta_{t_1},h\theta_{t_2})+\tilde{\mu}(\theta_{t_1},h\theta_{t_2}) \\
&-\tilde{\mu}(h\theta_{t_1},fg\theta_{t_2})+\tilde{\mu}(h\theta_{t_1},hd_\cP\theta_{t_2})+\tilde{\mu}(h\theta_{t_1},\theta_{t_2}).
\end{split}\end{align}

The first and third term are just $-\theta_{t,\epsilon_1}^\circ+\theta_{t,\epsilon_1}^{\Id}$, while the fourth and sixth are $-\theta_{t,\epsilon_2}^\circ+\theta_{t,\epsilon_2}^{\Id}$.

Also for a $t$ as above we have
\begin{equation}\label{eqn:three}
\theta_t d_{\cF\cP}=\tilde{\mu}(h\theta_{t_1}d_{\cF\cP},h\theta_{t_2})+\tilde{\mu}(h\theta_{t_1},h\theta_{t_2}d_{\cF\cP}).
\end{equation}
Now the second term of (\ref{eqn:two}) and the first term of (\ref{eqn:three}) are, by induction, equal to
$$
\sum_{\epsilon\in\mathbf{e}(t_1)}\theta_{t,\epsilon}^{\Id}-\theta_{t,\epsilon}^\circ,
$$
while the fifth term of (\ref{eqn:two}) and the second term of (\ref{eqn:three}) are equal to
$$
\sum_{\epsilon\in\mathbf{e}(t_2)}\theta_{t,\epsilon}^{\Id}-\theta_{t,\epsilon}^\circ
$$
and we are done.

If, say, $t_2$ is trivial we get
\begin{align*}\begin{split}
d_\cP \theta_t=&d_\cP \tilde{\mu}(h\theta_{t_1},\Id_\cP) \\
=&-\tilde{\mu}(d_\cP h\theta_{t_1},\Id_\cP)-\tilde{\mu}(h\theta_{t_1},d_\cP) \\
=&-\tilde{\mu}(fg\theta_{t_1},\Id_\cP)+\tilde{\mu}(hd_\cP \theta_{t_1},\Id_\cP)+\tilde{\mu}(\theta_{t_1},\Id_\cP)-\tilde{\mu}(h\theta_{t_1},d_\cP),\\
\end{split}\end{align*}
where the first and third term are $$-\theta_{t,\epsilon_1}^\circ+\theta_{t,\epsilon_1}^{\Id}$$ and by induction
$$
\tilde{\mu}(hd_\cP \theta_{t_1},\Id_\cP)+\left(\theta_t d_{\cF\cP}-\tilde{\mu}(h\theta_{t_1},d_\cP)\right)=\sum_{\epsilon\in\mathbf{e}(t_1)}\theta_{t,\epsilon}^{\Id}-\theta_{t,\epsilon}^\circ.
$$
\epro

\bex
Let $G$ be a three-vertex graph. As noted before, $T_G$ consists of two trees,
$$
T_G=\left\{
t_1=
\xygraph{
!{<0pt,0pt>;<10pt,0pt>:<0pt,10pt>::},
!{(0,0)}-!{(-2,2)},
!{(0,0)}-!{(2,2)},
!{(-1,1)}-!{(0,2)},
!{(0,0)}-!{(0,-1)},
},
t_2=
\xygraph{
!{<0pt,0pt>;<10pt,0pt>:<0pt,10pt>::},
!{(0,0)}-!{(-2,2)},
!{(0,0)}-!{(2,2)},
!{(1,1)}-!{(0,2)},
!{(0,0)}-!{(0,-1)},
}
\right\}.
$$
We denote each of the unique internal edges of $t_1$ and $t_2$ by $\epsilon$. Let $\times$ denote the differential $d_\cP$. A pictorial description of the proof in this special case is as follows:
\begin{align*}
d_\cP\theta_{t_1}=
\xygraph{
!{<0pt,0pt>;<10pt,0pt>:<0pt,10pt>::},
!{(0,0)}-!{(-2,2)},
!{(0,0)}-!{(2,2)},
!{(-1,1)}-!{(0,2)},
!{(0,0)}-!{(0,-1)},
!{(0,-0.5)}*{\times},
!{(-0.5,0.5)}*{\bullet}
}
=&
-\xygraph{
!{<0pt,0pt>;<10pt,0pt>:<0pt,10pt>::},
!{(0,0)}-!{(-2,2)},
!{(0,0)}-!{(2,2)},
!{(-1,1)}-!{(0,2)},
!{(0,0)}-!{(0,-1)},
!{(-0.3,0.3)}*{+},
!{(-0.7,0.7)}*{\bullet}
}
-
\xygraph{
!{<0pt,0pt>;<10pt,0pt>:<0pt,10pt>::},
!{(0,0)}-!{(-2,2)},
!{(0,0)}-!{(2,2)},
!{(-1,1)}-!{(0,2)},
!{(0,0)}-!{(0,-1)},
!{(1,1)}*{+},
!{(-0.5,0.5)}*{\bullet}
} \\
=&
\xygraph{
!{<0pt,0pt>;<10pt,0pt>:<0pt,10pt>::},
!{(0,0)}-!{(-2,2)},
!{(0,0)}-!{(2,2)},
!{(-1,1)}-!{(0,2)},
!{(0,0)}-!{(0,-1)},
!{(-0.3,0.3)}*{\bullet},
!{(-0.7,0.7)}*{+}
}
-
\xygraph{
!{<0pt,0pt>;<10pt,0pt>:<0pt,10pt>::},
!{(0,0)}-!{(-2,2)},
!{(0,0)}-!{(2,2)},
!{(-1,1)}-!{(0,2)},
!{(0,0)}-!{(0,-1)},
!{(-0.5,0.5)}*{\circ},
}
+
\xygraph{
!{<0pt,0pt>;<10pt,0pt>:<0pt,10pt>::},
!{(0,0)}-!{(-2,2)},
!{(0,0)}-!{(2,2)},
!{(-1,1)}-!{(0,2)},
!{(0,0)}-!{(0,-1)},
}
-
\xygraph{
!{<0pt,0pt>;<10pt,0pt>:<0pt,10pt>::},
!{(0,0)}-!{(-2,2)},
!{(0,0)}-!{(2,2)},
!{(-1,1)}-!{(0,2)},
!{(0,0)}-!{(0,-1)},
!{(1,1)}*{+},
!{(-0.5,0.5)}*{\bullet}
} \\
=&
\underbrace{
-
\xygraph{
!{<0pt,0pt>;<10pt,0pt>:<0pt,10pt>::},
!{(0,0)}-!{(-2,2)},
!{(0,0)}-!{(2,2)},
!{(-1,1)}-!{(0,2)},
!{(0,0)}-!{(0,-1)},
!{(-0.5,0.5)}*{\bullet},
!{(-1.5,1.5)}*{+}
}
-
\xygraph{
!{<0pt,0pt>;<10pt,0pt>:<0pt,10pt>::},
!{(0,0)}-!{(-2,2)},
!{(0,0)}-!{(2,2)},
!{(-1,1)}-!{(0,2)},
!{(0,0)}-!{(0,-1)},
!{(-0.5,0.5)}*{\bullet},
!{(-0.5,1.5)}*{+}
}
-
\xygraph{
!{<0pt,0pt>;<10pt,0pt>:<0pt,10pt>::},
!{(0,0)}-!{(-2,2)},
!{(0,0)}-!{(2,2)},
!{(-1,1)}-!{(0,2)},
!{(0,0)}-!{(0,-1)},
!{(1,1)}*{+},
!{(-0.5,0.5)}*{\bullet}
}}_{-\theta_{t_1}d_{\cF\cP}}
\underbrace{
-
\xygraph{
!{<0pt,0pt>;<10pt,0pt>:<0pt,10pt>::},
!{(0,0)}-!{(-2,2)},
!{(0,0)}-!{(2,2)},
!{(-1,1)}-!{(0,2)},
!{(0,0)}-!{(0,-1)},
!{(-0.5,0.5)}*{\circ},
}}_{-\theta_{t_1,\epsilon}^\circ}
+
\underbrace{
\xygraph{
!{<0pt,0pt>;<10pt,0pt>:<0pt,10pt>::},
!{(0,0)}-!{(-2,2)},
!{(0,0)}-!{(2,2)},
!{(-1,1)}-!{(0,2)},
!{(0,0)}-!{(0,-1)},
}}_{\theta_{t_1,\epsilon}^{\Id}},
\end{align*}
and in the same way
$$
d_\cP\theta_{t_2}=
\xygraph{
!{<0pt,0pt>;<10pt,0pt>:<0pt,10pt>::},
!{(0,0)}-!{(-2,2)},
!{(0,0)}-!{(2,2)},
!{(1,1)}-!{(0,2)},
!{(0,0)}-!{(0,-1)},
!{(0,-0.5)}*{\times},
!{(0.5,0.5)}*{\bullet}
}
=
\underbrace{
-
\xygraph{
!{<0pt,0pt>;<10pt,0pt>:<0pt,10pt>::},
!{(0,0)}-!{(-2,2)},
!{(0,0)}-!{(2,2)},
!{(1,1)}-!{(0,2)},
!{(0,0)}-!{(0,-1)},
!{(0.5,0.5)}*{\bullet},
!{(-1,1)}*{+}
}
-
\xygraph{
!{<0pt,0pt>;<10pt,0pt>:<0pt,10pt>::},
!{(0,0)}-!{(-2,2)},
!{(0,0)}-!{(2,2)},
!{(1,1)}-!{(0,2)},
!{(0,0)}-!{(0,-1)},
!{(0.5,0.5)}*{\bullet},
!{(0.5,1.5)}*{+}
}
-
\xygraph{
!{<0pt,0pt>;<10pt,0pt>:<0pt,10pt>::},
!{(0,0)}-!{(-2,2)},
!{(0,0)}-!{(2,2)},
!{(1,1)}-!{(0,2)},
!{(0,0)}-!{(0,-1)},
!{(1.5,1.5)}*{+},
!{(0.5,0.5)}*{\bullet}
}}_{-\theta_{t_2}d_{\cF\cP}}
\underbrace{
-
\xygraph{
!{<0pt,0pt>;<10pt,0pt>:<0pt,10pt>::},
!{(0,0)}-!{(-2,2)},
!{(0,0)}-!{(2,2)},
!{(1,1)}-!{(0,2)},
!{(0,0)}-!{(0,-1)},
!{(0.5,0.5)}*{\circ},
}}_{-\theta_{t_2,\epsilon}^\circ}
+
\underbrace{
\xygraph{
!{<0pt,0pt>;<10pt,0pt>:<0pt,10pt>::},
!{(0,0)}-!{(-2,2)},
!{(0,0)}-!{(2,2)},
!{(1,1)}-!{(0,2)},
!{(0,0)}-!{(0,-1)},
}}_{\theta_{t_2,\epsilon}^{\Id}}.
$$
As stated before, by associativity of $\mu$ we have $\theta_{t_1,\epsilon}^{\Id}+\theta_{t_2,\epsilon}^{\Id}=0$, and hence,
$$
d_\cP(\theta_{t_1}+\theta_{t_2})+(\theta_{t_1}+\theta_{t_2})d_{\cF\cP}=-\theta_{t_1,\epsilon}^\circ-\theta_{t_2,\epsilon}^\circ.
$$
\eex

\subsection{Induced morphism of sh properads}\hspace{1pt}

As in the classical theory on algebras we now show that the morphism $f:\cE\ra\cP$ extends to a morphism of sh properads.

\bde
We define $F_G:\decor{G}{\cE[1]}\ra\cP[1]$ to be $f$ on corollas, and for $|\mathbf{v}(G)|=k\geq 2$
$$
F_G=\sum_{t\in T_G}h\theta_t f^{\otimes k}.
$$
\ede

The family $F_G$ constitute a degree zero morphism $\bar\cF^c(\cE[1])\ra\cP[1]$, which extends to a well-defined morphism $F:\bar{\cF}^c(\cE[1])\ra B\cP$ of coproperads. Next we show that $F$ respects the differentials.

\bpr\label{prop:one}
In the same situation as in Theorem \ref{thm:one}, the morphism \\ $F:\bar{\cF}^c(\cE[1])\ra B\cP$ determined by $F_G$ is a morphism of dg coproperads.
\epr
\bpro
To prove that $F \partial_\cE=\partial_\cP F$, we examine the compositions with projection onto $\cP[1]$. If $(F \partial_\cE)_1=(\partial_\cP F)_1$, then the full morphisms agree since the equalities
\begin{align*}
\Delta F \partial_\cE&=(F,F)\Delta \partial_\cE=(F \partial_\cE,F)\Delta+(F,F \partial_\cE)\Delta\text{ and} \\
\Delta \partial_\cP F&=(\partial_\cP,\Id)\Delta F+(\Id,\partial_\cP)\Delta F=(\partial_\cP F,F)\Delta+(F,\partial_\cP F)\Delta
\end{align*}
show that the $k$-vertex components are determined by the $i$-vertex components for $i<k$, and are equal.

For $(F \partial_\cE)_1$ we have, on a decorated graph $G$ with $k$ vertices, that
\begin{align*}
(F \partial_\cE)_1&=\sum_{H\subset G}F_{G/H}\partial_H \\
&=\sum_{t\in T_G}fg\theta_t f^{\otimes k}+h\theta_t d_{\cF\cP}f^{\otimes k}+\sum_{\epsilon\in\mathbf{e}(t)}h\theta_{t,\epsilon}^\circ f^{\otimes k} \\
&=f\partial_G+h\theta_G d_{\cF\cP}f^{\otimes k}+\sum_{t\in T_G}\sum_{\epsilon\in\mathbf{e}(t)}h\theta_{t,\epsilon}^\circ f^{\otimes k}.
\end{align*}

For $(\partial_\cP F)_1$ we note that since $(\partial_\cP)_1$ is non-zero only on graphs with one or two vertices, it is enough to consider these parts of the image of $F$. Hence,
$$
(\partial_\cP F)_1=(\partial_\cP)_1 F_1+(\partial_\cP)_1 F_2=d_\cP F_1+\tilde{\mu}F_2.
$$
Now, on a decorated graph $G$ we have
$$
F_2=\sum_{G_1\circ G_2=G}(F_{G_1},F_{G_2}),
$$
so that, whith $k_i=|\mathbf{v}(G_i)|$, we get
\begin{align*}
(\partial_\cP F)_1&=\sum_{t\in T_G}d_\cP h\theta_t f^{\otimes k}+\sum_{\substack{G_1\circ G_2=G \\ t_1\in T_{G_1} \\ t_2\in T_{G_2}}}\tilde{\mu}(h\theta_{t_1}f^{\otimes k_1},h\theta_{t_2}f^{\otimes k_2}) \\
&=\sum_{t\in T_G}d_\cP h\theta_t f^{\otimes k}+\theta_t f^{\otimes k} \\
&=\sum_{t\in T_G}fg\theta_t f^{\otimes k}-h d_\cP \theta_t f^{\otimes k} \\
&=f \partial_G-h d_\cP\theta_G f^{\otimes k}.
\end{align*}
Hence we are done if we can prove
$$
-h d_\cP\theta_G f^{\otimes k}=h\theta_G d_{\cF\cP}f^{\otimes k}+\sum_{t\in T_G}\sum_{\epsilon\in\mathbf{e}(t)}h\theta_{t,\epsilon}^\circ f^{\otimes k},
$$
but this follows from (\ref{eqn:one}).
\epro

\subsection{Transfer for sh properads}\hspace{1pt}

The proofs of the results of the previous two sections may be modified to the situation where $\cP$ is an sh properad. In this Section we describe how this is done.

If $\cP$ is an sh properad, we have a collection of degree one morphisms
$$
\mu_G:\decor{G}{\cP[1]}\ra\cP[1]
$$
such that $\sum_{H\subset G}\mu_{G/H}\mu_H=0$. Contraction sequences for a $k$-vertex graph are now of length \emph{at most} $k-1$ and the set of equivalence classes of such sequences, $\hat{T}_G$, is a set of rooted trees with $k$ leaves, not necessarily binary. With the same definition of $\partial_\cE$, with $T_G$ replaced by $\hat{T}_G$, that is,
\begin{align*}
\theta_G=&\sum_{t\in\hat{T}_G}\theta_t, \\
\partial_G=&g\theta_G f^{\otimes k}\text{ and} \\
\partial_\cE=&\sum_{H\subset G}\partial_G,
\end{align*}
we have the following theorem.
\bth\label{thm:two}
Let $\cP$ be an sh properad and $\cE$ a $\Sigma$-bimodule, $g:\cP\ra\cE$, $f:\cE\ra\cP$ be morphisms of $\Sigma$-bimodules such that $fg-\Id_\cP=d_\cP h+h d_\cP$, where $h:\cP\ra\cP$ is a morphism of degree $-1$. Then the morphism $\partial_\cE$ is a codifferential on $\bar{\cF}(\cE[1])$.
\eth
\bpro[Sketch of proof]
The proof is similar to the proof of Theorem \ref{thm:one}. Lemma \ref{lem:three} is no longer true since $\mu$ is not associative. Hence, the question is how the $\theta_{t,\epsilon}^{\Id}$ cancel out. Obviously, for each tree $t$ in $\hat{T}_G$ and internal edge $\epsilon$ in $t$, the tree $t/\epsilon$ is also in $\hat{T}_G$, since if we may contract a subgraph in two steps then we can do it in one. If $t$ corresponds to a contraction sequence of length $l$, equation (\ref{eqn:two}) reads
\begin{align*}
d_\cP\tilde{\mu}_t&(h\theta_{t_1},\ldots,h\theta_{t_k}) \\
=&-\tilde{\mu}_{G_l}(d_\cP h\theta_{t_1},\ldots,h\theta_{t_l})-\cdots-\tilde{\mu}_G(h\theta_{t_1},\ldots,d_\cP h\theta_{t_l})-\sum_{(t',\epsilon)}\theta_{t',\epsilon}^{\Id}
\end{align*}
where the sum ranges over pairs $(t',\epsilon)$ with $t'\in\hat{T}_G$, $\epsilon$ an internal edge of $t'$ connected to the root vertex and $t'/\epsilon=t$. The terms $-\theta_{t',\epsilon}^{\Id}$ cancel the ones coming from $d_\cP\theta_{t'}$ when commuting $d_\cP$ past the $h$ decorating $\epsilon$.
\epro

\bex
Consider again a three-vertex graph $G$. Now,
$$
\hat{T}_G=\left\{
t_1=
\xygraph{
!{<0pt,0pt>;<10pt,0pt>:<0pt,10pt>::},
!{(0,0)}-!{(-2,2)},
!{(0,0)}-!{(2,2)},
!{(-1,1)}-!{(0,2)},
!{(0,0)}-!{(0,-1)},
},
t_2=
\xygraph{
!{<0pt,0pt>;<10pt,0pt>:<0pt,10pt>::},
!{(0,0)}-!{(-2,2)},
!{(0,0)}-!{(2,2)},
!{(1,1)}-!{(0,2)},
!{(0,0)}-!{(0,-1)},
},
t_3=
\xygraph{
!{<0pt,0pt>;<10pt,0pt>:<0pt,10pt>::},
!{(0,0)}-!{(-2,2)},
!{(0,0)}-!{(2,2)},
!{(0,0)}-!{(0,2)},
!{(0,0)}-!{(0,-1)},
}
\right\}.
$$
Also, $\theta_{t_3}=\mu_G$ and since $\cP$ is an sh properad,
$$
d_\cP\theta_{t_3}=
\xygraph{
!{<0pt,0pt>;<10pt,0pt>:<0pt,10pt>::},
!{(0,0)}-!{(-2,2)},
!{(0,0)}-!{(2,2)},
!{(0,0)}-!{(0,2)},
!{(0,0)}-!{(0,-1)},
!{(0,-0.5)}*{\times},
}
=
\underbrace{
-
\xygraph{
!{<0pt,0pt>;<10pt,0pt>:<0pt,10pt>::},
!{(0,0)}-!{(-2,2)},
!{(0,0)}-!{(2,2)},
!{(0,0)}-!{(0,2)},
!{(0,0)}-!{(0,-1)},
!{(-1,1)}*{+},
}
-
\xygraph{
!{<0pt,0pt>;<10pt,0pt>:<0pt,10pt>::},
!{(0,0)}-!{(-2,2)},
!{(0,0)}-!{(2,2)},
!{(0,0)}-!{(0,2)},
!{(0,0)}-!{(0,-1)},
!{(0,1)}*{\times},
}
-
\xygraph{
!{<0pt,0pt>;<10pt,0pt>:<0pt,10pt>::},
!{(0,0)}-!{(-2,2)},
!{(0,0)}-!{(2,2)},
!{(0,0)}-!{(0,2)},
!{(0,0)}-!{(0,-1)},
!{(1,1)}*{+},
}}_{-\theta_{t_3}d_{\cF\cP}}
\underbrace{
-
\xygraph{
!{<0pt,0pt>;<10pt,0pt>:<0pt,10pt>::},
!{(0,0)}-!{(-2,2)},
!{(0,0)}-!{(2,2)},
!{(-1,1)}-!{(0,2)},
!{(0,0)}-!{(0,-1)},
}
-
\xygraph{
!{<0pt,0pt>;<10pt,0pt>:<0pt,10pt>::},
!{(0,0)}-!{(-2,2)},
!{(0,0)}-!{(2,2)},
!{(1,1)}-!{(0,2)},
!{(0,0)}-!{(0,-1)},
}}_{-\theta_{t_1,\epsilon}^{\Id}-\theta_{t_2,\epsilon}^{\Id}}.
$$
We see that the last two terms are precisely the ones we need to cancel out the corresponding ones from $d_\cP\theta_{t_1}$ and $d_\cP\theta_{t_2}$, cf.~the example above. Hence,
$$
d_\cP(\theta_{t_1}+\theta_{t_2}+\theta_{t_3})+(\theta_{t_1}+\theta_{t_2}+\theta_{t_3})d_{\cF\cP}=-\theta_{t_1,\epsilon}^\circ-\theta_{t_2,\epsilon}^\circ.
$$
\eex

The proof of the proposition on the induced morphism of sh properads also remains valid. With the same definition of $F_G$, with $T_G$ replaced by $\hat{T}_G$, that is,
$$
F_G=\sum_{t\in\hat{T}_G}h\theta_t f^{\otimes k},
$$
the following proposition holds true.
\bpr\label{prop:two}
In the same situation as in Theorem \ref{thm:two}, the morphism \\ $F:\bar{\cF}^c(\cE[1])\ra B\cP$ determined by $F_G$ is a morphism of dg coproperads.
\epr
\bpro
The proof is the same as for Proposition \ref{prop:one}, we just note that
$$
(\partial_\cP F)_1=\sum_{l=1}^k(\partial_\cP)_1 F_l
$$
where
$$
(\partial_\cP)_1 F_1=\sum_{t\in\hat{T}_G}d_\cP h\theta_t f^{\otimes k}
$$
and
$$
\sum_{l=2}^k(\partial_\cP)_1 F_l=\sum_{t\in\hat{T}_G}\theta_t f^{\otimes k}.
$$
\epro

\bibliographystyle{amsalpha}
\bibliography{shprops}

\end{document}